\documentclass[11pt]{article}
\usepackage[letterpaper]{geometry}
\usepackage{amsmath,amsthm,amsfonts,amssymb}
\usepackage{enumerate,color}
\usepackage{graphicx}
\usepackage{subfigure}

\numberwithin{equation}{section}
\theoremstyle{plain}

\newtheorem{theorem}{Theorem}

\newtheorem{proposition}[theorem]{Proposition}

\begin{document}

\begin{center}
  \Large \bf Limit theorems for Markovian Hawkes processes with a large initial intensity
\end{center}

\author{}
\begin{center}
{Xuefeng
  Gao}\,\footnote{Department of Systems
    Engineering and Engineering Management, The Chinese University of Hong Kong, Shatin, N.T. Hong Kong;
    xfgao@se.cuhk.edu.hk},
  Lingjiong Zhu\,\footnote{Corresponding Author. Department of Mathematics, Florida State University, 1017 Academic Way, Tallahassee, FL-32306, United States of America; zhu@math.fsu.edu.
  }
\end{center}

\begin{center}
 \today
\end{center}

\begin{abstract}
Hawkes process is a simple point process that is self-exciting
and has clustering effect. The intensity of this point process depends on its
entire past history. It has wide applications in finance, neuroscience, social networks, criminology, seismology, and many
other fields. In this paper, we study the linear Hawkes process with an exponential exciting function in the asymptotic regime where the initial intensity of the Hawkes process is large.  We derive limit theorems under this asymptotic regime as well as the regime when both the initial intensity and the time are large.
\end{abstract}

\section{Introduction}

Let $N$ be a simple point process on $\mathbb{R}$ and let $\mathcal{F}^{-\infty}_{t}:=\sigma(N(C),C\in\mathcal{B}(\mathbb{R}), C\subset(-\infty,t])$ be
an increasing family of $\sigma$-algebras. Any nonnegative $\mathcal{F}^{-\infty}_{t}$-progressively measurable process $\lambda_{t}$ with
\begin{equation*}
\mathbb{E}\left[N(a,b]|\mathcal{F}^{-\infty}_{a}\right]=\mathbb{E}\left[\int_{a}^{b}\lambda_{s}ds\big|\mathcal{F}^{-\infty}_{a}\right] \quad \text{almost surely,}
\end{equation*}
 for all intervals $(a,b]$ is called an $\mathcal{F}^{-\infty}_{t}$-intensity of $N$. We use the notation $N_{t}:=N(0,t]$ to denote the number of
points in the interval $(0,t]$.

We consider $N$ to be a linear Hawkes process, which is a simple point process admitting an $\mathcal{F}^{-\infty}_{t}$-intensity
\begin{equation}
\lambda_{t}:=\mu + \int_{-\infty}^{t-}h(t-s)N(ds) =  \mu + \sum_{\tau_i<t}h(t-\tau_i),\label{dynamics}
\end{equation}
where $\mu \ge 0$ is the base intensity, $\tau_i$ are the occurrence times of the points before time $t$, and
$h(\cdot):\mathbb{R}^{+}\rightarrow\mathbb{R}^{+}$ is the exciting function encoding the influence of past events on the intensity and
we always assume that $\Vert h\Vert_{L^{1}}=\int_{0}^{\infty}h(t)dt<\infty$ and $h$ is locally bounded.

The linear Hawkes process was first introduced by A.G. Hawkes
in 1971 \cite{Hawkes}. It naturally generalizes the Poisson process and it
captures both the self--exciting property and the clustering effect. In addition, Hawkes process is a very versatile model which is amenable to statistical analysis. These explain why it has wide applications in neuroscience, genome analysis,
criminology, social networks, seismology, insurance, finance and many other fields.
For a list of references, we refer to \cite{ZhuThesis}.

Throughout this paper, we assume an exponential exciting function $h(t):=\alpha e^{-\beta t}$ where $\alpha,\beta>0$. That is, we restrict ourselves
to the linear Markovian Hawkes process\footnote{Note that, the Hawkes process is in general non--Markovian. When the exciting function is given by a sum of exponential functions, the Hawkes process $(\lambda, N)$ is Markovian \cite{BacrySurvey}.}. To see the Markov property, we define
\begin{equation} \label{eq:Z-def}
Z_{t}:=\int_{-\infty}^{t}\alpha e^{-\beta(t-s)}N(ds) = Z_0 \cdot e^{-\beta t} + \int_{0}^{t}\alpha e^{-\beta(t-s)}N(ds) .
\end{equation}
Then, the process $Z$ is Markovian and satisfies the dynamics:
\begin{equation}\label{eq:markov-Z}
dZ_{t}=-\beta Z_{t}dt+\alpha dN_{t},
\end{equation}
where $N$ is a Hawkes process with intensity $\lambda_t = \mu + Z_{t-}$ at time $t$. In addition, the pair $(Z, N)$ is also Markovian. We also assume $Z_0=Z_{0-}$, i.e., there is no jump at time zero.

We consider an asymptotic regime where $Z_{0}=n$, and $n \in \mathbb{R}^{+}$ is sent to infinity.
We derive functional law of large numbers and functional central limit theorems for linear Markovian Hawkes processes in this asymptotic regime as well as the regime when both the initial intensity and the time are large.

Our main results (Theorem~1--4) are mainly oriented to develop approximations for the transient behavior of Hawkes processes with large initial intensity $\lambda_0$. Note that $\lambda_0 =\mu +Z_{0}$, so when $Z_{0}$ is large, we have $\lambda_0$ large. In practice, this means the intensity of arriving events observed at time zero is considerably larger than the base intensity $\mu$. Our results could be potentially useful for applications where the Hawkes process is a relevant model. As an example, consider stock trading in limit order markets. Arrivals of trades or orders are usually clustered in time and can be modeled by Markovian Hawkes
processes, see, e.g., Bowsher \cite{Bowsher}, Da Fonseca and Zaatour \cite{DaFonseca2014}, Hewlett \cite{Hewlett2006}.
If the current trading intensity of a stock is high (a given large value), then our asymptotic analysis can help approximate the volume duration, i.e., the waiting time until a predetermined buy (or sell) volume is traded\footnote{Assuming trade size is a constant given by the average trade size over a certain time window.}. Such a volume duration is useful in measuring the (time) costs of liquidity in the trading process, see, e.g., Gourieroux et al. \cite{Gourieroux1999}, Hautsch \cite{Hautsch}.
As an other example, consider portfolio credit risk. Markovian Hawkes
processes have been proposed to capture corporate default clustering in a portfolio of interacting firms, see, e.g., Errias et al. \cite{Errais}. Our asymptotic analysis can help approximate the number of corporate defaults occurring over a given time interval
when the initial intensity of defaults is high such as during a financial crisis.

To the best of our knowledge, this is the first paper to study the large initial intensity asymptotics
in the context of the Hawkes process. For simplicity, the discussions in our paper are restricted to the case when the exciting function is a single exponential function $h(t)=\alpha e^{-\beta t}$. Indeed, our results can be extended to the case when the exciting function is a sum of exponentials so that the Hawkes process is still Markovian. Specifically, when $h(t)=\sum_{i=1}^{d}\alpha_{i}e^{-\beta_{i}t}$, where $\beta_{i}>0$, $\alpha_{i}\in\mathbb{R}$
such that $h(t)>0$ for every $t$, $(Z_{t}^{1},\ldots,Z_{t}^{d})$ is a $d$-dimensional Markov process,
where $Z_{t}^{i}=\int_{-\infty}^{t}\alpha_{i}e^{-\beta_{i}(t-s)}dN_{s}$, see e.g. \cite{ZhuI}.
The results in this paper for large $Z_0$ asymptotics can be extended to a Hawkes process with large initial values of $(Z_{0}^{1},\ldots,Z_{0}^{d})$
as $Z_{0}^{i}$ goes to infinity at the same order of speed.
We do not study here the non--Markovian case which may require a different approach since knowing the value of the initial intensity is not
sufficient to describe the future evolution of a non--Markovian Hawkes process.
On the other hand, any reasonable exciting function $h(t)$ can be approximated by a linear combination of exponential functions,
see e.g. \cite{ZhuI}. In this respect, the Markovian setting in this paper is not too restrictive. From the application point of view, Markovian Hawkes processes are the most widely used due to the tractability of the theoretical analysis
as well as the simulations and calibrations. See, e.g., Bacry et al. \cite{BacrySurvey} and the references therein.

\textbf{Related Literature.}
We now explain the difference between our work and the existing literature.
Note that, almost all the existing literature on limit theorems
for Hawkes processes are for large time asymptotics, e.g., the functional law of large numbers
and functional central limit theorems in Bacry et al. \cite{Bacry}, the large deviations principle in Bordenave
and Torrisi \cite{Bordenave} and the moderate deviation principle in Zhu \cite{ZhuMDP} for linear Hawkes processes. These large--time
 limit theorems
hold under the so--called subcritical regime, that is when $\Vert h\Vert_{L^{1}}<1$.
The critical and supercritical regimes are when $\Vert h\Vert_{L^{1}}=1$ and $\Vert h\Vert_{L^{1}}>1$ respectively. See, e.g., \cite{ZhuThesis, TorrisiI} for limit theorems for nonlinear Hawkes processes under these different regimes. The limit theorems for the nearly critical, or so--called nearly unstable case where $\Vert h\Vert_{L^{1}}$ is close to one was studied
in Jaisson and Rosenbaum \cite{Jaisson,JRII}.

Other than the large time asymptotics, the large dimensional asymptotics of Hawkes processes have been studied recently,
see e.g. Delattre et al. \cite{Delattre}, that is, the asymptotics
for the multivariate Hawkes process where the number of dimension goes to infinity.
The mean-field limit was obtained in \cite{Delattre} for different regimes.
See also Chevallier \cite{Chevallier}, Hodara and L\"{o}cherbach \cite{HL}, Delattre and Fournier \cite{DF}. Our work complements all these studies by studying the asymptotic behavior of a linear Markovian Hawkes process with a large initial intensity.

\textbf{Organization of this paper.} The rest of the paper is organized as follows. In Section 2, we state our main results. In Section 3, we present the proofs of the main results. The proofs of auxiliary results are collected in the appendix.

\section{Main results}
In this section we state our main results. Note that processes $Z$ and $N$ both depend on the initial condition $Z_{0}=n$ and
one can use $Z^{(n)},N^{(n)}$ to emphasize the dependence on $Z_{0}=n$.
For simplicity of notations, we use $Z_{t},N_{t}$ to denote $Z_{t}^{(n)},N_{t}^{(n)}$
throughout this paper. Write $D[a,b]$ as the space of c\`{a}dl\`{a}g processes on $[a,b] \subset[0, \infty)$ that are equipped with Skorohod $J_1$ topology (see e.g., Billingsley \cite{Billingsley}).

\subsection{Limit theorems with large $Z_0$}
In this section, we present limit theorems including functional law of large numbers (FLLN) and functional central limit theorem (FCLT) for
processes $Z$ and $N$ in the asymptotic regime where $Z_{0}=n \rightarrow \infty$.

We first state the FLLN for the processes $Z$ and $N$.

\begin{theorem} \label{thm:LLN}
Fix any $T>0$. As $n\rightarrow\infty$, we have

\begin{eqnarray}
\sup_{0\leq t\leq T}\left|\frac{Z_{t}}{n}-e^{(\alpha-\beta)t}\right|\rightarrow 0,
\qquad
\text{almost surely}, \label{eq:LLN-Z}\\
\sup_{0\leq t\leq T}\left|\frac{N_{t}}{n}-\psi(t)\right|\rightarrow 0,
\qquad
\text{almost surely}, \label{eq:LLN-N}
\end{eqnarray}
where
\begin{equation}\label{eq:psi}
\psi(t):= \begin{cases}
\frac{e^{(\alpha-\beta)t}-1}{\alpha-\beta},
 & \alpha \ne \beta, \\
t, & \alpha = \beta.
\end{cases}
\end{equation}

\end{theorem}

We next state the FCLT for the processes $Z$ and $N$. Recall that a Gaussian process is called centered if its mean function is identically zero.

\begin{theorem} \label{thm:CLT}
For any $T>0$, as $n\rightarrow\infty$, the sequence of processes
\begin{align}
\left\{ \frac{Z_{t}-ne^{(\alpha-\beta)t}}{\sqrt{n}}: t \in [0,T] \right\} \label{eq:FCLT-Z}
\end{align}
converges in distribution to a centered Gaussian process $G$ on $D[0,T]$,
where the covariance function of $G$ is given as follows: for $0\le s \le t$,
\begin{equation}\label{eq:covG}
\text{Cov}(G_s,G_t)= \begin{cases}
\frac{\alpha^{2}}{\alpha-\beta}
(e^{(\alpha-\beta)(t+s)}-e^{(\alpha-\beta)t}),  & \alpha \ne \beta, \\
\alpha^2 s, & \alpha = \beta.
\end{cases}
\end{equation}
Furthermore, the sequence of re--normalized Hawkes processes
\begin{align}
\left\{\frac{N_{t}-n \cdot \psi(t)}{\sqrt{n}}: t \in [0,T] \right\} \label{eq:FCLT-N}
\end{align}
converges in distribution to a centered Gaussian process $H$ on $D[0,T]$, where $H$ is given by
\begin{eqnarray} \label{eq:H}
H_t:=\frac{G_t}{\alpha}+\frac{\beta}{\alpha}\int_{0}^{t}G_s ds, \quad t \in [0,T].
\end{eqnarray}
\end{theorem}

We next discuss two properties of the limiting Gaussian processes $G$ and $H$.

First,  when $\alpha = \beta$, one readily finds from \eqref{eq:covG} that the limiting Gaussian process $G$ is actually a Brownian motion with drift zero and variance $\alpha^2$, and its sample paths are (almost surely) continuous. When $\alpha \ne \beta$, one can also readily verify  from \eqref{eq:covG} and Kolmogrov's continuity criterion that the limiting Gaussian process $G$ with covariance \eqref{eq:covG} has a version with continuous sample paths on the bounded time interval $[0,T]$. Hence the paths of the process $H$ in \eqref{eq:H} inherit such continuity property.

Second,
the Gaussian process $G$ is a Markov process and $H$ is not. This can be easily verified from the criteria that a centered Gaussian process $\Upsilon$ with covariance function $\Gamma(s,t):=\mathbb{E}[\Upsilon_s \Upsilon_t]$ is Markovian if and only if
\begin{eqnarray*}
\Gamma(s,u) \Gamma(t, t) = \Gamma(s,t) \Gamma(t,u).
\end{eqnarray*}
for every $0 \le s < t < u$. See, e.g., Revuz and Yor \cite[p.86]{Revuz}.

To illustrate the usefulness of our limit theorems, let us consider the first passage time problem mentioned in the introduction. For a given volume $K$, we are interested in the waiting
time $\tau_K$ defined by $\tau_K :=\inf\{t>0: N_t \ge K\}$ where by convention the infimum of an empty set is $+\infty$. When $Z_0=n$ is large, we can approximate the distribution of $\tau_K$ using Theorem~\ref{thm:CLT} as follows: for each $0<t<\infty$,
\begin{eqnarray} \label{eq:prob}
\mathbb{P}(\tau_K \le t ) = \mathbb{P}(N_t \ge K) \approx \mathbb{P}(n \psi(t) + \sqrt{n}H_t \ge K) = 1- \Phi\left(\frac{K- n \psi(t)}{\sqrt{n \cdot Var(H_t)}} \right),
\end{eqnarray}
where $\psi(\cdot)$ is given in \eqref{eq:psi}, $Var(H_t)$ is the variance of $H_t$ and $\Phi(\cdot)$ is the cumulative distribution function of a standard normal random variable. The normal approximation in \eqref{eq:prob} could work well if $K- n \psi(t)$ is on the order of $\sqrt{n}$ for large $n$. In such a case, the probability in \eqref{eq:prob} can be readily computed after one fits a Markovian Hawkes processes to the data (observations of occurrence times of events before time zero) to estimate the parameters $\alpha$ and $\beta$, and then obtain $Z_0=n$ from \eqref{eq:Z-def} using the occurrence times of past events. See, e.g., Ozaki \cite{Ozaki}, Daley and
Vere--Jones \cite[Chapter~7]{Daley} for estimation methods and further details.

\subsection{Limit theorems with large $Z_0$ and large time}
In this section, we present limit theorems for
processes $Z$ and $N$ in the asymptotic regime where both $Z_{0}=n$ and the time go to infinity.
Such limit theorems could provide insights on the `macroscopic' behavior of the Hawkes processes with large initial intensity.

When the time is sent to infinity, Hawkes processes behave differently depending on the value of
$\Vert h\Vert_{L^{1}}$ (see, e.g., Zhu \cite{ZhuThesis}).
In our case,  the exciting function is exponential: $h(t)= \alpha e^{-\beta t}$. So we have the following different regimes:
(a) Super--critical: $\alpha> \beta$;
(b) Sub--critical: $\alpha < \beta$;
 (c) Critical: $\alpha= \beta$;
(d) Nearly--critical: $\alpha \approx \beta$ (the precise definition will be given in Theorem~\ref{thm:rescaled-CLT}).
We study each case separately.

We first state the FLLN for the process $Z$ and the point process $N$ when $Z_0$ and the time are both large.
\begin{theorem} \label{thm:rescale-LLN}

\begin{itemize}
\item [(i)] (Super--Critical Case) Assume that $\alpha>\beta>0$ and let $\tau_{n}=\frac{\log n}{\alpha-\beta}$.
For any $T>0$, as $n\rightarrow\infty$ we have
\begin{eqnarray}
\sup_{0\leq s\leq T}\left|\frac{Z_{s \tau_{n}}}{n^{1+s}}-1\right|\rightarrow 0, \qquad \text{almost surely,} \label{eq:Zlargetime}\\
\sup_{0\leq s\leq T}\left|\frac{N_{s \tau_{n}}}{n^{1+s}}-\frac{1}{\alpha-\beta}+\frac{1}{(\alpha-\beta)n^{s}}\right|\rightarrow 0, \qquad \text{almost surely.} \label{eq:Nlargetime}
\end{eqnarray}

\item [(ii)] (Sub--Critical Case) Assume that $\beta>\alpha>0$ and let $t_{n}=\frac{\log n}{\beta-\alpha}$.
Then, for any $0<T<1$, as $n\rightarrow\infty$ we have
\begin{eqnarray}
\sup_{0\leq s\leq T}\left|\frac{Z_{st_{n}}}{n^{1-s}}-1\right|\rightarrow 0, \qquad \text{in probability,} \label{eq:Zlargetime-sub}\\
\sup_{0\leq s\leq T}\left|\frac{N_{st_{n}}}{n}-\frac{1}{\beta-\alpha}+\frac{1}{(\beta-\alpha)n^{s}}\right|\rightarrow 0, \label{eq:Nlargetime-sub}
\qquad
\text{in probability}.
\end{eqnarray}
We have almost surely convergence for \eqref{eq:Zlargetime-sub} and \eqref{eq:Nlargetime-sub} if $0<T<\frac{1}{2}$.
\end{itemize}
\end{theorem}
A few remarks are in order. First, we choose to speed up the time by a factor of $\log n$ and scale the space accordingly to ensure the convergence
to a non--trivial limit.
Such a scaling simplifies our presentation and it is natural in view of \eqref{eq:LLN-Z} and \eqref{eq:LLN-N} in Theorem~\ref{thm:LLN}.
Second, in the sub--critical case, we restrict $T <1$.
One can easily show that $\mathbb{E}[Z_{st_{n}}]=n^{1-s}$ (see Proposition~\ref{prop:moments-noncritial}), which
implies that $\mathbb{E}[Z_{st_{n}}]=1$ if $s=1$.
Therefore, the process $Z$ starts with a large number $Z_0=n$, and after $t_n$ unit of time, it drops to one on average ($Z$ is always positive).
So it makes sense to consider the law of large numbers
for $Z_{st_{n}}/{n^{1-s}}$ for $0\leq s\leq T$ where $0<T<1$. Due to our proof technique, the almost surely convergence here is restricted to $T< \frac{1}{2}$.
Third, we do not state the FLLN in the critical or nearly--critical case since the limit process is trivial under appropriate time and space scalings, as readily seen from the result below.

We next state the FCLT for the processes $Z$ and $N$ when $Z_0=n$ is large and the time is speeded up in an appropriate way.
\begin{theorem} \label{thm:rescaled-CLT}
(i)(Critical and Nearly Critical Cases) Fix $\beta>0$ and $\gamma \in \mathbb{R}$. For $Z_0 =n$, define $\alpha_n = \beta + \frac{\gamma}{n}$, which is positive for all large $n$.
Assume that
\begin{equation*}  \label{eq:Z-critical}
dZ_{t}= - \beta Z_{t}dt+\alpha_n dN_{t},
\end{equation*}
where the point process $N$ has intensity $\mu +Z_{t-}$ at time $t$.
Then as $n \rightarrow \infty$, we have the sequence of processes $\left\{\frac{Z_{tn}}{n}: t \in [0,T] \right\}$
converges in distribution to the process $X$ on $D[0,T]$, \text{where $X$} satisfies
\begin{align}
dX_{t}=(\beta \mu + \gamma X_{t})dt+\beta \sqrt{X_{t}}dB_{t},\qquad X_{0}=1, \label{eq:CIR}
\end{align}
where $B$ is a standard Brownian motion. In addition, we have that 
the sequence of re--normalized Hawkes processes $\left\{\frac{N_{tn}}{n^2}: t \in [0,T] \right\}$
converges in distribution toward the process \[ \left\{\int_{0}^{t}X_{s}ds: t \in [0,T] \right\}\] on $D[0,T]$.

(ii)(Super--Critical Case) Assume that $\alpha>\beta>0$ and let $\tau_{n}=\frac{\log n}{\alpha-\beta}$. For any $0<t<T<\frac{1}{2}$, as $n \rightarrow \infty$, we have the sequence of processes
\begin{equation}\label{eq:Zlargetime-super-clt}
\left\{\frac{Z_{s \tau_{n}}-n^{1+s}}{\sqrt{n^{1+2s}}}, \quad s \in [t,T]\right\},
\end{equation}
converges in distribution to the process $Y$ on $D[t,T]$, where $Y_s \equiv \xi$ for $s \in [t,T]$ and $\xi$ is a normal random variable with mean $0$ and variance $\frac{\alpha^{2}}{\alpha-\beta}$.
Moreover, we have
\begin{equation}\label{eq:Nlargetime-super-clt}
\left\{\frac{N_{s \tau_{n}}-\frac{n^{1+s}-n}{\alpha-\beta}}{n^{\frac{1}{2}+s}}, \quad s \in [t,T]\right\},
\end{equation}
converges in distribution to the process $\frac{Y}{\alpha - \beta}$ on $D[t,T]$.

(iii) (Sub--Critical Case) Assume that $\beta>\alpha>0$ and let $t_{n}=\frac{\log n}{\beta-\alpha}$. For any $0<T<1$, we have all the finite dimensional distributions of the processes
\begin{equation} \label{eq:sub-largetime}
\left\{\frac{Z_{st_{n}}-n^{1-s}}{\sqrt{n^{1-s}}}, \quad s \in [0,T]\right\},
\end{equation}
converges in distribution to the corresponding finite dimensional distributions of a centered Gaussian process $R$ with covariance function $\text{Cov}(R_u, R_v)=\frac{\alpha^{2}}{\beta-\alpha}$ if $u=v>0$ and zero otherwise.
\end{theorem}

We provide some further discussions on limit processes and the convergence mode in Theorem~\ref{thm:rescaled-CLT}.

First, since $\beta \mu \ge 0$
and $X_0=1$, there is a unique strong solution $X \ge 0$ to the stochastic differential equation \eqref{eq:CIR}, see, e.g., \cite[Chapter IX]{Revuz}. The process $X$ is known as a square--root diffusion or a Cox--Ingersoll--Ross (CIR) process, which is widely used in modeling interest rate and stochastic volatility in financial mathematics. In addition, the properties of $X$ and the time integral of $X$ are well studied,
see e.g.  \cite{going2003survey}. In the critical case where $\gamma =0$, the process $\frac{4}{\beta^2}X$  with $X$ defined by \eqref{eq:CIR} is also known as the square of a $\frac{4 \mu}{\beta}$--dimensional Bessel process.

In the nearly--critical case when $\gamma<0$, we have $0<\alpha_n < \beta$ for each large $n$ and $\alpha_n/\beta$ approaches one as $n \rightarrow \infty$. Such Hawkes processes
are called nearly unstable in Jaisson and Rosenbaum \cite{Jaisson}. In line with the results in \cite{Jaisson}, we obtain, in this nearly unstable case, mean--reverting CIR process in the limit for the process $Z$ where $Z_0$ and the time are both sent to infinity. We also remark that since we consider a different asymptotic regime compared with \cite{Jaisson}, we obtain a different initial condition for the limiting process $X$ in \eqref{eq:CIR}. In addition, our result holds for general $\gamma \in \mathbb{R}$.

Second, in the super--critical case, the weak limit of the sequence of processes in \eqref{eq:Zlargetime-super-clt} is a constant normal random variable $\xi$.
Note that this process--level convergence is restricted to $D[t, T]$ for $T>t>0$, and such weak convergence can not be extended to the space $D[0,T]$, which is readily seen after noting that $Z_0=n$.

Third, in the sub--critical case, the centered Gaussian process $R$, with covariance function $\text{Cov}(R_u, R_v)=\frac{\alpha^{2}}{\beta-\alpha}$ if $u=v$ and zero otherwise, is known to exist. But such a Gaussian process does not have a measurable version and hence a c\`{a}dl\`{a}g version, see, e.g. Revuz and Yor \cite[p.37]{Revuz}.
Therefore, we only have the convergence of finite dimensional distributions for the sequence of processes in \eqref{eq:sub-largetime}, but we do not have process--level convergence. In other words, the sequence of processes in \eqref{eq:sub-largetime} is not tight. Hence we also do not have or state weak convergence result for the Hawkes process $N$ (when both $Z_0$ and the time are large) in this sub--critical case.

Finally, we remark that unlike the critical and nearly--critical case cases,
we choose to speed up the time by a factor of $\log n$ in the super and sub--critical cases.  As in Theorem~\ref{thm:rescale-LLN}, such time scalings simplify the notations in the paper.

\section{Proofs of Main Results}
In this section we gather the proofs of our main results Theorem~1--4. For notational simplicity, in all the proofs we use $C>0$ as a generic constant which may vary line from line.
The constant $C$ may depend on $\alpha, \beta, \gamma$ and $T$, but it is
independent of $n$.

Our proof strategy is to first show the functional law of large numbers and functional central limit theorems for processes $Z$ and $N$
when the base intensity $\mu$ is zero, and then extend the proofs to the case when $\mu>0$.
Such a strategy relies critically on the
observation described in the next section.

\subsection{Decomposition of Hawkes processes} \label{sec:mu-relation}
This section presents a decomposition result for linear Hawkes processes.

The linear Hawkes process has the well--known immigration birth representation, see, e.g., \cite{HawkesII}
\footnote{In \cite{HawkesII}, they assume
$\Vert h\Vert_{L^{1}}<1$ because they are considering the stationary regime.
In our setting, we only consider a fixed time interval $[0,t]$ given a finite
initial intensity, and the immigration birth representation still holds without
the constraint $\Vert h\Vert_{L^{1}}<1$.}.
That is, the immigrant arrives according to a homogeneous
Poisson process with constant rate $\mu$. Each immigrant would produce children and the number
of children has Poisson distribution with parameter $\Vert h\Vert_{L^{1}}=\frac{\alpha}{\beta}$. Conditional on the number of
the children of an immigrant, the time that a child was born
has probability density function $\frac{h(t)}{\Vert h\Vert_{L^{1}}}=\beta e^{-\beta t}$. Each child would produce children
according to the same laws independent of other children. All the immigrants produce children independently.
The number of points of a linear Hawkes process on a time interval $[0,t]$ equals
the total number of immigrants and the descendants on the interval $[0,t]$.

Recall that we are interested in a Hawkes process $N$ with intensity $\mu+Z_{t-}$ at time $t$ and $Z_{0}=n$, and
\begin{equation*}
Z_{t}=ne^{-\beta t}+\int_{0}^{t}\alpha e^{-\beta(t-s)}dN_{s}.
\end{equation*}
By the immigration-birth representation, we can decompose linear Hawkes process as
\begin{eqnarray} \label{eq:decomp-N}
N_{t}=N_{t}^{(0)}+N_{t}^{(1)}, \quad \text{$t \ge 0$},
\end{eqnarray}
where $N_{t}^{(0)}$ is the number of points of the immigrants that arrive according
to an inhomogeneous Poisson process with rate $ne^{-\beta t}$ and all the descendants
on the interval $[0,t]$
and $N_{t}^{(1)}$ is the number of points of the immigrants that arrive according
to a homogeneous Poisson process with rate $\mu$ and all the descendants on the interval $[0,t]$.
Therefore,
$N^{(0)}$ is a simple point process with intensity $Z^{(0)}$, where
\begin{equation*}
Z_{t}^{(0)}=ne^{-\beta t}+\int_{0}^{t}\alpha e^{-\beta(t-s)}dN_{s}^{(0)},
\end{equation*}
and in the differential form,
\begin{equation*}
dZ_{t}^{(0)}=-\beta Z_{t}^{(0)}dt+\alpha dN_{t}^{(0)}, \quad Z_{0}^{(0)}=n.
\end{equation*}
$N^{(1)}$ is a simple point process with intensity $\lambda^{(1)}$ where
\begin{equation*}
\lambda^{(1)}_{t}:=\mu+Z_{t}^{(1)}=\mu+\int_{0}^{t}\alpha e^{-\beta(t-s)}dN_{s}^{(1)},
\end{equation*}
where
\begin{equation*}
Z_{t}^{(1)}=\int_{0}^{t}\alpha e^{-\beta(t-s)}dN_{s}^{(1)},
\end{equation*}
and in the differential form,
\begin{equation*} \label{eq:Z1}
dZ_{t}^{(1)}=-\beta Z_{t}^{(1)}dt+\alpha dN_{t}^{(1)},\quad Z_{0}^{(1)}=0.
\end{equation*}
In addition, the two point processes $N^{(0)}$ and $N^{(1)}$ are independent of each other. As a result, we also have
\begin{eqnarray} \label{eq:decomp-Z}
Z_{t}=Z_{t}^{(0)}+Z_{t}^{(1)}, \quad \text{$t \ge 0$},
\end{eqnarray}
and the processes $Z^{(0)}$ and $Z^{(1)}$ are also independent of each other. An alternative
way to see the validity of decompositions \eqref{eq:decomp-N} and \eqref{eq:decomp-Z} is via the Poisson embedding technique often used in the theory of point processes \cite{Bremaud}.

Now we illustrate the high level idea of our proofs of limit theorems. The first step is to study the Hawkes process $N$ with $\mu=0$. In the decomposition $N=N^{(0)}+N^{(1)}$,
we note that $N^{(1)}$ is a linear Hawkes process which is empty on $(-\infty, 0]$. Thus when $\mu=0$, $N_t^{(1)} \equiv 0$ for any $t \ge 0$. So we have $N= N^{(0)}$ and $Z=Z^{(0)}$ in this case. Once we have established limit theorems for $Z^{(0)}$ and $N^{(0)}$, we move to the second step: consider Hawkes process with $\mu>0$. In view of the decompositions in \eqref{eq:decomp-N} and \eqref{eq:decomp-Z}, it suffices to have limit theorems for $N^{(1)}$ and $Z^{(1)}$ when $n \rightarrow \infty$. This is relatively straightforward: when we consider a finite time horizon, then $N^{(1)}$ and $Z^{(1)}$ do not contribute
in the scaling limits as they are independent of $Z_0 =Z_0^{(0)}= n$ which goes to infinity; When we consider the time also goes to infinity, the functional law of large numbers and central limit theorems of $N^{(1)}$
have already been well studied in the literature. So
combining the limit theorems we obtained for $N^{(0)}$ and $Z^{(0)}$ in the first step,
we have the limit theorems for processes $N$ and $Z$.

\subsection{Preliminaries}
This section presents preliminaries for proving the main results. Throughout this section,
$N$ is a simple point process and its intensity at time $t$ is given by $Z_{t-}$
where $dZ_{t}=-\beta Z_{t}dt+\alpha dN_{t}$. That is, we consider $\mu=0$.

We have the following immediate observation. We can express
the jump process $N_{t}=N(0,t]$ in terms of the intensity process $Z_{t}$ in the following way,
\begin{equation} \label{eq:N-Z}
N_{t}=\frac{Z_{t}-Z_{0}}{\alpha}+\frac{\beta}{\alpha}\int_{0}^{t}Z_{s}ds.
\end{equation}
This relation is very useful: once we prove limit theorems for process $Z$, the identity \eqref{eq:N-Z} allows us to prove limit theorems for process $N$ in a direct way.

Define for $t \ge 0,$
\begin{eqnarray*}
M_{t}:=N_{t}-\int_{0}^{t}Z_{s}ds.
\end{eqnarray*}
Then $M$ is a martingale with respect to the natural filtration, and the predictable quadratic variation $\langle M\rangle_{t}$
of this martingale is given by $\int_{0}^{t}Z_{s}ds$, see, e.g. \cite{Lipster}. It is readily seen that
\begin{align*}
d(Z_{t}-ne^{(\alpha-\beta)t})
&=-\beta Z_{t}dt+\alpha dN_{t}-(\alpha-\beta)ne^{(\alpha-\beta)t}dt
\\
&=(\alpha-\beta)(Z_{t}-ne^{(\alpha-\beta)t})dt+\alpha dM_{t}.
\nonumber
\end{align*}
Therefore, we have
\begin{equation}\label{MartingaleRep}
Z_{t}-ne^{(\alpha-\beta)t}=e^{(\alpha-\beta)t}\alpha\int_{0}^{t}e^{-(\alpha-\beta)s}dM_{s}.
\end{equation}

We next summarize two results on the moments of $Z_t$ and related inequalities. The proofs are given in the appendix.
\begin{proposition} \label{prop:moments-noncritial}
Suppose $\alpha,\beta > 0$ and $\mu=0$. We have
\begin{eqnarray*}
\mathbb{E}[Z_{t}|Z_{0}]=Z_{0}e^{(\alpha-\beta)t}.
\end{eqnarray*}
If $\alpha = \beta,$ we have
\begin{eqnarray*}
\mathbb{E}[Z_{t}^{2}|Z_0] &=&
Z_{0}^{2}+\alpha^{2}Z_{0}t, \\
\mathbb{E}[Z_{t}^{3}|Z_0]&=&Z_{0}^{3}+3\alpha^{2}Z_{0}^{2}t+\frac{3}{2}\alpha^{4}Z_{0}t^{2}+\alpha^{3}Z_{0}t.
\end{eqnarray*}
If $\alpha \ne \beta$, we have
\begin{eqnarray*}
\mathbb{E}[Z_{t}^{2}|Z_0] &=& Z_{0}^{2}e^{2(\alpha-\beta)t}
+\frac{\alpha^{2}Z_{0}}{\alpha-\beta}(e^{2(\alpha-\beta)t}-e^{(\alpha-\beta)t}),\\
\mathbb{E}[Z_{t}^{3} |Z_0]
&=&\left(Z_{0}^{3}+\frac{3\alpha^{2}Z_{0}^{2}}{\alpha-\beta}
+\frac{\alpha^{3}Z_{0}}{2(\alpha-\beta)} + \frac{3\alpha^{4}Z_{0}}{2(\alpha-\beta)^{2}}\right)
e^{3(\alpha-\beta)t}
\\
&&\qquad
-\left(\frac{3\alpha^{2}Z_{0}^{2}}{\alpha-\beta}+\frac{3\alpha^{4}Z_{0}}{(\alpha-\beta)^{2}}\right)
e^{2(\alpha-\beta)t}
-\left(\frac{\alpha^{3}Z_{0}}{2(\alpha-\beta)}-\frac{3\alpha^{4}Z_{0}}{2(\alpha-\beta)^{2}}\right)
e^{(\alpha-\beta)t}.
\nonumber
\end{eqnarray*}
\end{proposition}

\begin{proposition}\label{prop:ineq}
Given $Z_0 =n$ sufficiently large, we have
\begin{eqnarray}
\sup_{0 \le t \le T}\mathbb{E}\left[\left(Z_{t}-ne^{(\alpha-\beta)t}\right)^{4}\right]
&\leq & C n^{2},\label{UpperBoundI} \\
\sup_{\delta \le t \le T}\mathbb{E}\left[\left(\int_{t-\delta}^{t}e^{-(\alpha-\beta)s}dM_{s}\right)^{4}\right]
&\leq& C n^{2}\delta^2 , \quad \text{for $\delta \in [0,t]$.}\label{eq:mart-upperbound}
\end{eqnarray}
In addition, when $\alpha \ne \beta$ we have for any $t \ge 0$
\begin{align} \label{ineq:mart}
&\mathbb{E}\left[\left(\int_{0}^{t}e^{-(\alpha-\beta)s}dM_{s}\right)^{4}\right]
\\
&\leq Ct\left[\frac{n^{2}}{2(\alpha-\beta)}(1-e^{-2(\alpha-\beta)t})
+\frac{\alpha^{2}n}{2(\alpha-\beta)^{2}}(1-e^{-2(\alpha-\beta)t})
-\frac{\alpha^{2}n}{3(\alpha-\beta)^{2}}(1-e^{-3(\alpha-\beta)t})\right].
\nonumber
\end{align}
Here $C$ is a constant that may depend on $\alpha, \beta, \gamma$ and $T$, but it is
independent of $n$ and $t$.
\end{proposition}

\subsubsection{Moment generating function of $Z_t$}
In this section we present a result on moment generating function of $Z_t$,
which can be found in Errais et al. \cite{Errais} for example. We include it here since it is a critical tool in establishing Gaussian limit in the functional central limit theorems. Indeed, in \cite{Errais}, they only discussed the Laplace transforms.
But for the purpose of showing the convergence of distributions, we need convergence
of moment generating functions in a neighborhood of zero.

As $Z$ is an affine process, we can infer that (see, e.g., \cite{Errais, Keller2015}) for $\theta \in i \mathbb{R}$, the function $u(z,t):=\mathbb{E}[e^{-\theta Z_{t}}|Z_{0}=z]$ is given by
 $u(z, t) = e^{A(t)z}$, where
\begin{align}
&A'(t)=-\beta A(t)+e^{A(t)\alpha}-1, \label{eq:A}\\
& A(0)=-\theta. \label{eq:A0}
\end{align}
To show the real exponential moment of $Z_t$ exists in a neighborhood of zero, we use \cite[Theorem~2.14]{Keller2015} and study the ODE for $A$ as given above. One can readily verify that for $\theta \in \mathbb{R}$ and $t>0$, the ODE system \eqref{eq:A}--\eqref{eq:A0} has an unique real--valued solution that starts at $-\theta < \theta_c(t)$ and exists up to time $t$ (i.e., the solution does not blow up on $[0,t]$)\footnote{Note that in \cite{Keller2015}, 
they use the notion of the minimal solution since for the most general affine process,
with the presence of L\'{e}vy measures, the solutions of the ODE system reaching or starting at the boundary
of the domain may not be unique. In our case, $-\beta A+e^{A\alpha}-1$ is locally Lipschitz continuous
in $A$ on $\mathbb{R}$, and the local solution of the ODE system is unique.}, 
where $\theta_c(t) \in (0, \infty)$ is defined by
\begin{equation}\label{thetacEqn}
\theta_{c}(t):=\sup\left\{\theta>0:\int_{\theta}^{\infty}\frac{dA}{-\beta A+e^{\alpha A}-1}=t\right\},
\end{equation}
which depends on $t, \alpha$ and $\beta$. Hence, we deduce from \cite[Theorem~2.14]{Keller2015} that for $-\theta\in (-\infty, \theta_{c}(t))$,
\begin{eqnarray}\label{eq:u}
u(z,t)=\mathbb{E}[e^{-\theta Z_{t}}|Z_{0}=z] = e^{A(t, -\theta)z} < \infty.
\end{eqnarray}
For convenience, here and in the following we write $A(t, -\theta)$ instead of $A(t)$ to emphasize that $A$ takes value $-\theta$ at time 0.
We remark that in general the function $A$ is not explicit in the sense that there is no closed-form solution to the differential equation system given by \eqref{eq:A} and \eqref{eq:A0}.
The key ideas in establishing our limit theorems include exploiting perturbation theory (see, e.g., \cite{murdock1999perturbations}) of differential equations and Gronwall's inequality to obtain estimates of ODE solutions in order to show the convergence of certain moment generating functions.

\subsection{Proof of Theorem~\ref{thm:LLN}}
We prove Theorem~\ref{thm:LLN} in this section. We first consider the case $\mu=0$, and then extend the proof to the case $\mu>0$ using the observation in Section~\ref{sec:mu-relation}.

\subsubsection{Proof of Theorem~\ref{thm:LLN} when $\mu=0$}

\begin{proof} [Proof of \eqref{eq:LLN-Z}]
When $\mu=0$, let us recall from \eqref{MartingaleRep} that
\begin{equation*}
Z_{t}-ne^{(\alpha-\beta)t}=e^{(\alpha-\beta)t}\alpha\int_{0}^{t}e^{-(\alpha-\beta)s}dM_{s}.
\end{equation*}
Together with Doob's martingale inequality, we have for any $\epsilon>0$,
\begin{eqnarray*}
\mathbb{P}\left(\sup_{0\leq t\leq T}\left|\frac{Z_{t}}{n}-e^{(\alpha-\beta)t}\right|\geq\epsilon\right)
&\leq&\mathbb{P}\left(\alpha\left(1+e^{|\alpha-\beta|T}\right)\sup_{0\leq t\leq T}
\left|\int_{0}^{t}e^{-(\alpha-\beta)s}dM_{s}\right|\geq n\epsilon\right)
\nonumber
\\
&\leq&\frac{\alpha^{4}(1+e^{|\alpha-\beta|T})^{4}}{n^{4}\epsilon^{4}}
\mathbb{E}\left[\left(\int_{0}^{T}e^{-(\alpha-\beta)s}dM_{s}\right)^{4}\right]. \nonumber
\end{eqnarray*}
On combining with inequality \eqref{ineq:mart} we obtain
\begin{align*}
\mathbb{P}\left(\sup_{0\leq t\leq T}\left|\frac{Z_{t}}{n}-e^{(\alpha-\beta)t}\right|\geq\epsilon\right) \leq\frac{Cn^2 + C n}{n^{4}\epsilon^{4}}. \nonumber
\end{align*}
Since $\sum_{n=1}^{\infty} \frac{Cn^2 + C n}{n^{4}\epsilon^{4}}$ is finite, the result \eqref{eq:LLN-Z} follows from the Borel-Cantelli lemma.
\end{proof}

\begin{proof}[Proof of \eqref{eq:LLN-N}]
By \eqref{eq:LLN-Z}, we have
\begin{equation*}
\sup_{0\leq t\leq T}\left|\frac{\int_{0}^{t}Z_{s}ds}{n}-\psi(t)\right|\rightarrow 0,
\qquad
\text{almost surely as $n\rightarrow\infty$},
\end{equation*}
where $\psi(t)=\int_{0}^{t}e^{(\alpha-\beta)s}ds$. Since $M_{t}=N_{t}-\int_{0}^{t}Z_{s}ds$, thus in order to show \eqref{eq:LLN-N}, it suffices to show
that
\begin{eqnarray} \label{eq:mt0}
\sup_{0\leq t\leq T}\frac{|M_{t}|}{n}\rightarrow 0,
\qquad
\text{almost surely as $n\rightarrow\infty$}.
\end{eqnarray}
Similar to the proof of \eqref{eq:LLN-Z}, we can apply Doob's martingale inequality and \eqref{ineq:mart} to
show that for any $\epsilon>0$
\begin{align}
&\mathbb{P}\left(\sup_{0\leq t\leq T}\frac{|M_{t}|}{n} \ge \epsilon\right)
\leq\frac{C}{n^{4}\epsilon^{4}}
\mathbb{E}M_T^4
\nonumber \leq \frac{1}{n^4 \epsilon^4} \cdot C (n^2 + n ).
\nonumber
\end{align}
Thus \eqref{eq:mt0} follows from Borel-Cantelli lemma.

\subsubsection{Proof of Theorem~\ref{thm:LLN} when $\mu>0$}
As described in Section~\ref{sec:mu-relation}, when $\mu>0$, we can decompose $Z_{t}=Z_{t}^{(0)}+Z_{t}^{(1)}$ where $Z^{(0)}$ and $Z^{(1)}$ are independent. We have established in the previous section that
\begin{equation*}
\sup_{0\leq t\leq T}\left|\frac{Z_{t}^{(0)}}{n}-e^{(\alpha-\beta)t}\right|\rightarrow 0, \quad \text{almost surely.}
\end{equation*}
In addition, note
that $Z_{t}^{(1)}$ is independent of $Z_{0}^{(0)}=n$ and hence $\sup_{0\leq t\leq T}Z_{t}^{(1)}/n \rightarrow 0$ almost surely as $n \rightarrow \infty$.
Now \eqref{eq:LLN-Z} immediately follows.

Similarly, we have $N_{t}=N_{t}^{(0)}+N_{t}^{(1)}$ when $\mu>0$. Since $N_{t}^{(1)}$ is independent of the parameter $n$, we obtain $N_{T}^{(1)}/n \rightarrow 0$ almost surely as $n \rightarrow \infty$.
Thus \eqref{eq:LLN-N} follows.

\subsection{Proof of Theorem~2}
We prove Theorem~2 in this section.

\subsubsection{FCLT for $Z$ when $\mu=0$}

In this section we prove the weak convergence of re--normalized processes of $Z$
in \eqref{eq:FCLT-Z} on $D[0,T]$ for the case $\mu=0.$
For notational simplicity, we define for $Z_0=n$ and each $t \ge 0$
\begin{equation}\label{eq:diff-Z}
\tilde Z_{t} := \frac {Z_{t}-ne^{(\alpha-\beta)t} }{\sqrt{n}}.
\end{equation}
Our approach is to apply Theorem 13.5 in Billingsley \cite{Billingsley}. In particular, we verify the three conditions: (recall that $G$ is a centered Gaussian process with covariance \eqref{eq:covG})
\begin{itemize}
\item [(a)] $G_{t}-G_{t-\Delta}$ converge to zero in distribution as $\Delta \rightarrow 0$.
\item [(b)] Finite--dimensional distributions of $\tilde Z$ converge to those of $G$.
\item [(c)]  For $0 \le t -\delta \le t + \delta \le T$, there exists a constant $C$ independent of $n$ such that
\begin{equation*}
\mathbb{E}\left[(\tilde{Z}_{t+\delta}-\tilde{Z}_{t})^{2}(\tilde{Z}_{t}-\tilde{Z}_{t-\delta})^{2}\right]
\leq C\delta^{2}.
\end{equation*}
\end{itemize}

We first prove (a). When $\alpha = \beta$, we know $G$ is a Brownian motion with mean zero and variance $\alpha^2$. Hence
$\mathbb{E}[(G_{t}-G_{t-\Delta})^{2}]= \alpha^2 \Delta$ which goes to
zero as $\Delta\rightarrow 0$. Hence, by Chebychev's inequality, we conclude
that $G_{t}-G_{t-\Delta}\rightarrow 0$ in distribution
as $\Delta\rightarrow 0$. When $\alpha \ne \beta$, we find from \eqref{eq:covG} that $\mathbb{E}[(G_{t}-G_{t-\Delta})^{4}] \le C \Delta^2$. Hence
$G_{t}-G_{t-\Delta}\rightarrow 0$ in distribution
as $\Delta\rightarrow 0$ by Chebychev's inequality.

We next prove (b).
We first show that for each fixed $t \ge 0,$ the sequence of random variables $\{\tilde Z_{t}: n \ge 1\}$ defined in \eqref{eq:diff-Z} converges in distribution as $n \rightarrow \infty$. To this end, we study the moment generating function of $\tilde Z_{t}$. Fix $\theta \in \mathbb{R}$. It is immediate from \eqref{eq:u} that, for any sufficiently large $n$,
\begin{equation} \label{eq:mgf1}
\mathbb{E}\left[e^{-\theta\frac{Z_{t}-ne^{(\alpha-\beta)t}}{\sqrt{n}}}\bigg|Z_{0}=n\right]
=\exp \left({ A\left(t, -\frac{\theta}{\sqrt{n}}\right) \cdot n+\sqrt{n}\theta e^{(\alpha-\beta)t}} \right),
\end{equation}
To show this sequence of moment generating functions converges when $n\rightarrow \infty,$ we rely on the expansion of $A$ which satisfies the ordinary differential equation \eqref{eq:A} with initial condition $-\frac{\theta}{\sqrt{n}}$. Note that for $n$ large, the quantity $-\frac{\theta}{\sqrt{n}}$ is small. As the ODE solution $A$ depends smoothly on its initial value,
we introduce the following expansion: 
\begin{eqnarray} \label{eq:A-approx}
A\left(t, -\frac{\theta}{\sqrt{n}}\right)
= f_0(t) + f_1(t) \cdot \left(-\frac{\theta}{\sqrt{n}}\right) + f_2(t) \cdot \frac{\theta^2}{n} + O(n^{-\frac{3}{2}}),
\end{eqnarray}
where $O(\epsilon)$ is a term bounded by $C_t \cdot \epsilon$ for a positive constant $C_t$ and $\epsilon$ small enough. Due to the smooth dependence of the ODE for $A$ on its initial value, we infer that $f_0, f_1, f_2$ and the constant in the big O notation are all uniformly bounded for $t \in [0, T]$
 (see, e.g., \cite{murdock1999perturbations} for background on perturbation theory and asymptotic expansions for differential equations).
Next
we use the differential equation \eqref{eq:A} to determine the unknown functions  $f_0, f_1$ and $f_2$. First, it is obvious that $f_0(t) = A(t,0) \equiv 0$, i.e., the solution to \eqref{eq:A} is zero when we have zero initial condition. As we have $A(t, -\frac{\theta}{\sqrt{n}}) = O(\frac{1}{\sqrt{n}})$ uniformly for $t \in [0, T]$, we can deduce that
\begin{equation*}
e^{\alpha A\left(t, -\frac{\theta}{\sqrt{n}}\right) } - \beta A\left(t, -\frac{\theta}{\sqrt{n}}\right) -1
= (\alpha -\beta) A\left(t, -\frac{\theta}{\sqrt{n}}\right)
+\frac{1}{2}\alpha^{2}A\left(t, -\frac{\theta}{\sqrt{n}}\right)^{2} +O(n^{-\frac{3}{2}}).
\end{equation*}
Together with \eqref{eq:A-approx}, the differential equation \eqref{eq:A} becomes
\begin{eqnarray*}
\lefteqn{f_1'(t) \cdot \left(-\frac{\theta}{\sqrt{n}}\right) + f_2'(t) \cdot \frac{\theta^2}{n} + O(n^{-\frac{3}{2}})} \\
&=(\alpha -\beta) f_1(t) \cdot \left(-\frac{\theta}{\sqrt{n}}\right) 
+\left[\frac{1}{2}\alpha^{2}f_1(t)^2 + (\alpha -\beta) f_2(t)\right] \cdot \frac{\theta^2}{n} + O(n^{-\frac{3}{2}}).
\end{eqnarray*}
Hence equating the coefficients depending on the power of $-\frac{\theta}{\sqrt{n}}$,  we have
\begin{eqnarray*}
f_1'(t) &=& (\alpha -\beta) f_1(t),\\
f_2'(t) &=& \frac{1}{2}\alpha^{2}f_1(t)^2 + (\alpha -\beta) f_2(t).
\end{eqnarray*}
For the initial conditions, we can deduce from the equality $A(0, -\frac{\theta}{\sqrt{n}})= -\frac{\theta}{\sqrt{n}}$ that
$f_1(0)  = 1$, and $f_2(0)=0$.
Solving the above two differential equations we obtain
\begin{eqnarray}
f_1(t) &=& e^{(\alpha -\beta)t}, \label{eq:f1} \nonumber\\
f_2(t) &=& \begin{cases}
\frac{1}{2}\frac{\alpha^{2}}{\alpha-\beta}
(e^{2(\alpha-\beta)t}-e^{(\alpha-\beta)t}), & \text{if} \quad \alpha \ne \beta, \\
\frac{1}{2} \alpha^2 t, & \text{if} \quad \alpha = \beta.
\end{cases}
\label{eq:f2}
\end{eqnarray}
Thus we obtain
\begin{equation} \label{eq:A-approx2}
\quad A\left(t, -\frac{\theta}{\sqrt{n}}\right) =  e^{(\alpha -\beta)t} \cdot \left(-\frac{\theta}{\sqrt{n}}\right) + f_2(t) \cdot \frac{\theta^2}{n} + O\left(n^{-\frac{3}{2}}\right),
\end{equation}
which implies that as $n \rightarrow \infty$
\begin{equation*}
\left({ A\left(t, -\frac{\theta}{\sqrt{n}}\right) \cdot n+\sqrt{n}\theta e^{(\alpha-\beta)t}} \right)
\rightarrow \frac{\theta^{2}}{2} \cdot 2f_2(t).
\end{equation*}
Therefore, the sequences of the moment generating functions in \eqref{eq:mgf1} converges when $n\rightarrow \infty,$ and we obtain for fixed $t>0$, as $n \rightarrow \infty,$
\begin{equation*}
\frac{Z_{t}-ne^{(\alpha-\beta)t}}{\sqrt{n}}
\rightarrow G_{t}, \quad \text{in distribution},
\end{equation*}
where $G_{t}$ is a random variable following normal distribution with mean zero and variance $2 f_2(t)$ where $f_2$ is given in \eqref{eq:f2}.

Now to prove (b), it remains to show the finite dimensional distributions of $\tilde Z$ in \eqref{eq:diff-Z} converge for dimensions of two and higher. We use mathematical induction and rely on the Markov property of $Z$. Given $Z_0=n$, suppose for arbitrary $0 < t_1< \ldots <t_k$,
\begin{equation*}
\left(\tilde Z_{t_1}, \ldots, \tilde Z_{t_k} \right) \rightarrow \left(G_{t_1}, \ldots, G_{t_k}\right), \quad \text{in distribution.}
\end{equation*}
Here we suppose $\left(G_{t_1}, \ldots, G_{t_k}\right)$ follows $k-$variate normal distribution with mean zero and $k\times k$ covariance matrix $[\sigma_{ij}]_{k \times k}$ where for $1 \le i \le j \le k$
\begin{equation} \label{eq:covar}
\sigma_{ij} := \text{Cov}(G_{t_i},G_{t_j})= \begin{cases}
\frac{\alpha^{2}}{\alpha-\beta}
(e^{(\alpha-\beta)(t_i+t_j)}-e^{(\alpha-\beta)t_j}), & \alpha \ne \beta, \\
\alpha^{2} t_i, & \alpha = \beta ,
\end{cases}
\end{equation}
and $\sigma_{ij} = \sigma_{ji}$ if $i>j.$ Hence for $\theta_1, \ldots, \theta_k \in \mathbb{R}$, we have as $n \rightarrow \infty$
\begin{equation} \label{eq:mgf-k}
\mathbb{E}\left[e^{- \sum_{i=1}^{k}\theta_i \tilde Z_{t_i}}\big|Z_{0}=n\right]
\rightarrow \exp \left(\frac{1}{2} \sum_{i=1}^{k} \sum_{j=1}^{k} \theta_i \theta_j \sigma_{ij} \right).
\end{equation}
To show for arbitrary $0 < t_1< \ldots <t_k < t_{k+1}$,
\begin{equation*}
\left(\tilde Z_{t_1}, \ldots, \tilde Z_{t_{k+1}} \right) \rightarrow \left(G_{t_1}, \ldots, G_{t_{k+1}}\right), \quad \text{in distribution,}
\end{equation*}
where $\left(G_{t_1}, \ldots, G_{t_{k+1}}\right)$ follows $(k+1)-$variate normal distribution with mean zero and some covariance matrix consistent with \eqref{eq:covar},
we show the sequence of moment generating functions converges as $n \rightarrow \infty$. By tower property of conditional expectation and Markov property of $Z$, for any sufficiently large $n$, we have
\begin{equation} \label{eq:mgf-k1}
\mathbb{E}\left[e^{- \sum_{i=1}^{k+1}\theta_i \tilde Z_{t_i}}\bigg|Z_{0}=n\right] =
\mathbb{E}\left[e^{- \sum_{i=1}^{k}\theta_i \tilde Z_{t_i}} \cdot
\mathbb{E} \left[e^{- \theta_{k+1} \tilde Z_{t_{k+1}}}| Z_{t_{k}}\right] \bigg|Z_{0}=n\right].
\end{equation}
Similar for \eqref{eq:mgf1}, we first deduce from \eqref{eq:u} that
\begin{eqnarray*}
\lefteqn{\mathbb{E} \left[e^{- \theta_{k+1} \tilde Z_{t_{k+1}}}| Z_{t_{k}} \right] }\\
&=& \exp \left(A\left(t_{k+1} - t_k, -\frac{\theta_{k+1}}{\sqrt{n}}\right) Z_{t_{k}} + {\theta_{k+1}\sqrt{n}e^{(\alpha-\beta)t_{k+1}}}\right)\\
&=& \exp \left(A\left(t_{k+1} - t_k, -\frac{\theta_{k+1}}{\sqrt{n}}\right) \cdot ( \sqrt{n} \tilde Z_{t_{k}} + n e^{(\alpha - \beta) t_k} ) + {\theta_{k+1}\sqrt{n}e^{(\alpha-\beta)t_{k+1}}}\right),\\
\end{eqnarray*}
where the last equality follows from \eqref{eq:diff-Z}.
With this, we deduce from \eqref{eq:mgf-k1} that
\begin{eqnarray} \label{eq:mgf-k1-fin}
\mathbb{E}\left[e^{- \sum_{i=1}^{k+1}\theta_i \tilde Z_{t_i}}\big|Z_{0}=n\right] =
\mathbb{E}\left[e^{- \sum_{i=1}^{k} \hat \theta_i^{(n)} \tilde Z_{t_i}} \bigg|Z_{0}=n\right] \cdot \exp( \Gamma^n  ),
\end{eqnarray}
where $\hat \theta_i^{(n)}=\theta_i$ for $i=1, \ldots, k-1$, $\hat \theta_k^{(n)} = \theta_k - \sqrt{n} A(t_{k+1} - t_k, -\frac{\theta_{k+1}}{\sqrt{n}})$ and
\begin{eqnarray*}
\Gamma^n := A\left(t_{k+1} - t_k, -\frac{\theta_{k+1}}{\sqrt{n}}\right)  \cdot n e^{(\alpha-\beta)t_k} + \theta_{k+1} \sqrt{n} e^{(\alpha-\beta)t_{k+1}}.
\end{eqnarray*}
>From the expansion of $A$ at \eqref{eq:A-approx2} we infer that
\begin{eqnarray*}
\lim_{n \rightarrow \infty} \hat \theta_k^{(n)}  &=& \theta_k + \theta_{k+1} e^{(\alpha-\beta)(t_{k+1} - t_k)}, \\
\lim_{n \rightarrow \infty} \Gamma^n &=& \theta_{k+1}^2  e^{(\alpha-\beta)t_k} \cdot f_2 (t_{k+1} - t_k),
\end{eqnarray*}where $f_2$ is the function given in \eqref{eq:f2}.
In conjunction with \eqref{eq:mgf-k} and \eqref{eq:mgf-k1-fin}, and simplifying the resulting expression we obtain as $n \rightarrow \infty$,
\begin{equation*}
\mathbb{E}\left[e^{- \sum_{i=1}^{k+1}\theta_i \tilde Z_{t_i}}\big|Z_{0}=n\right]
\rightarrow \exp \left(\frac{1}{2} \sum_{i=1}^{k+1} \sum_{j=1}^{k+1} \theta_i \theta_j \sigma_{ij} \right),
\end{equation*}
where
\begin{eqnarray*}
\sigma_{(k+1)(k+1)} &=& 2 f_2 (t_{k+1}),\\
\sigma_{i(k+1)} &=& \begin{cases}
\frac{\alpha^{2}}{\alpha-\beta}
(e^{(\alpha-\beta)(t_i+t_{k+1})}-e^{(\alpha-\beta)t_{k+1}}), & \alpha \ne \beta, \\
\alpha^2 t_i, & \alpha = \beta,
\end{cases}
\end{eqnarray*}
and $\sigma_{(k+1)i} = \sigma_{i(k+1)}$ for $i=1, \ldots, k$.
Hence we have established that the finite--dimensional distributions of $\tilde Z$ in \eqref{eq:diff-Z} converges to that of a centered Gaussian process $G$ with covariance function is given by: for $0 \le s \le t$,
\begin{equation*}
\text{Cov}(G_{t},G_{s})= \begin{cases}
\frac{\alpha^{2}}{\alpha-\beta}
(e^{(\alpha-\beta)(t+s)}-e^{(\alpha-\beta)t})  & \alpha \ne \beta, \\
\alpha^2 s, & \alpha = \beta.
\end{cases}
\end{equation*}

Finally we prove (c).
We study the cases $\alpha \ne \beta$ and $\alpha = \beta$ separately.
We start with the case $\alpha \ne \beta$.
Given Proposition~\ref{prop:moments-noncritial}, it is straightforward to compute that
\begin{align} \label{eq:1}
&\mathbb{E}\left[\left(Z_{t+\delta}-ne^{(\alpha-\beta)(t+\delta)}-Z_{t}
+ne^{(\alpha-\beta)t}\right)^{2}\big|Z_{t}\right]
\\
&=\mathbb{E}[Z_{t+\delta}^{2}|Z_{t}]
-2(ne^{(\alpha-\beta)(t+\delta)}-ne^{(\alpha-\beta)t}+Z_{t})\mathbb{E}[Z_{t+\delta}|Z_{t}]
\nonumber
\\
&\qquad
+(ne^{(\alpha-\beta)(t+\delta)}-ne^{(\alpha-\beta)t}+Z_{t})^{2}
\nonumber
\\
&=(e^{(\alpha-\beta)\delta}-1)^{2}(Z_{t}-ne^{(\alpha-\beta)t})^{2}
+(e^{(\alpha-\beta)\delta}-1)\frac{\alpha^{2}e^{(\alpha-\beta)\delta}}{\alpha-\beta}Z_{t}.
\nonumber
\end{align}
In addition, by the Markovian property of $Z$ we obtain
\begin{align*}
&\mathbb{E}\left[(\tilde{Z}_{t+\delta}-\tilde{Z}_{t})^{2}(\tilde{Z}_{t}-\tilde{Z}_{t-\delta})^{2}\right]
\\
&= \mathbb{E}\left[ (\tilde{Z}_{t}-\tilde{Z}_{t-\delta})^{2} \cdot \mathbb{E}\left[(\tilde{Z}_{t+\delta}-\tilde{Z}_{t})^{2}\big| Z_t\right] \right].
\end{align*}
On combining the above equation with \eqref{eq:1} we obtain
\begin{align}\label{MidStep}
&\mathbb{E}\left[(\tilde{Z}_{t+\delta}-\tilde{Z}_{t})^{2}(\tilde{Z}_{t}-\tilde{Z}_{t-\delta})^{2}\right]
\\
&=\frac{1}{n^{2}}(e^{(\alpha-\beta)\delta}-1)^{2}\mathbb{E}\left[(Z_{t}-ne^{(\alpha-\beta)t})^{2}
\left(Z_{t}-ne^{(\alpha-\beta)t}-Z_{t-\delta}
+ne^{(\alpha-\beta)(t-\delta)}\right)^{2}\right]
\nonumber
\\
&\qquad
+\frac{1}{n^{2}}(e^{(\alpha-\beta)\delta}-1)\frac{\alpha^{2}e^{(\alpha-\beta)\delta}}{\alpha-\beta}
\mathbb{E}\left[Z_{t}\left(Z_{t}-ne^{(\alpha-\beta)t}-Z_{t-\delta}
+ne^{(\alpha-\beta)(t-\delta)}\right)^{2}\right].
\nonumber
\end{align}

Next we derive upper bounds for the terms in \eqref{MidStep}.
First, let us bound the first term on the right hand side of \eqref{MidStep}. Direct computation yields
\begin{align*}\label{eq:4}
&\mathbb{E}\left[(Z_{t}-ne^{(\alpha-\beta)t})^{2}
\left(Z_{t}-ne^{(\alpha-\beta)t}-Z_{t-\delta}
+ne^{(\alpha-\beta)(t-\delta)}\right)^{2}\right]
\\
&\leq 2\mathbb{E}\left[(Z_{t}-ne^{(\alpha-\beta)t})^{2}
\left[\left(Z_{t}-ne^{(\alpha-\beta)t}\right)^{2}+\left(Z_{t-\delta}
-ne^{(\alpha-\beta)(t-\delta)}\right)^{2}\right]\right]
\nonumber
\\
&\leq
2\mathbb{E}\left[\left(Z_{t}-ne^{(\alpha-\beta)t}\right)^{4}\right]
+2\left[\mathbb{E}\left(Z_{t}-ne^{(\alpha-\beta)t}\right)^{4}\right]^{1/2}
\left[\mathbb{E}\left(Z_{t-\delta}-ne^{(\alpha-\beta)(t-\delta)}\right)^{4}\right]^{1/2}
\nonumber
\\
&\leq 4Cn^{2},
\nonumber
\end{align*}
where the second last inequality follows from Cauchy-Schwarz inequality and the last inequality is due to Proposition~\ref{prop:ineq}. On combining with the fact that
\begin{eqnarray} \label{eq:6}
| e^{(\alpha -\beta) \delta} - 1 |  \le C \delta \quad \text{for} \quad \delta \in [0, T],
\end{eqnarray}
we infer that the first term on the right hand side of \eqref{MidStep} is upper bounded by $C \delta^2$ for some constant $C.$

We next proceed to bound the second term on the right hand side of \eqref{MidStep}.
By \eqref{MartingaleRep}, we have
\begin{align}
&Z_{t}-ne^{(\alpha-\beta)t}-(Z_{t-\delta}-ne^{(\alpha-\beta)(t-\delta)}) \nonumber
\\
&=e^{(\alpha-\beta)t}\alpha\int_{0}^{t}e^{-(\alpha-\beta)s}dM_{s}
-e^{(\alpha-\beta)(t-\delta)}\alpha\int_{0}^{t-\delta}e^{-(\alpha-\beta)s}dM_{s}
\nonumber
\\
&=e^{(\alpha-\beta)t}\alpha\int_{t-\delta}^{t}e^{-(\alpha-\beta)s}dM_{s}
+(e^{(\alpha-\beta)\delta}-1)e^{(\alpha-\beta)(t-\delta)}\alpha\int_{0}^{t-\delta}e^{-(\alpha-\beta)s}dM_{s}.
\nonumber
\end{align}
Therefore, by Cauchy-Schwarz inequality we find
\begin{align*}
&\mathbb{E}\left[Z_{t}\left(Z_{t}-ne^{(\alpha-\beta)t}-Z_{t-\delta}
+ne^{(\alpha-\beta)(t-\delta)}\right)^{2}\right]
\\
&\leq[\mathbb{E} Z_{t}^{2}]^{1/2}
\left[\mathbb{E}\left(Z_{t}-ne^{(\alpha-\beta)t}-Z_{t-\delta}
+ne^{(\alpha-\beta)(t-\delta)}\right)^{4}\right]^{1/2}
\nonumber
\\
&\leq[\mathbb{E} Z_{t}^{2}]^{1/2}
\bigg[4 \mathbb{E}\left(e^{(\alpha-\beta)t}\alpha\int_{t-\delta}^{t}e^{-(\alpha-\beta)s}dM_{s}\right)^{4}
\nonumber
\\
&\qquad\qquad\qquad
+4  \mathbb{E} \left((e^{(\alpha-\beta)\delta}-1)e^{(\alpha-\beta)(t-\delta)}\alpha\int_{0}^{t-\delta}e^{-(\alpha-\beta)s}dM_{s}\right)^{4}\bigg]^{1/2}
\nonumber
\\
&=2[\mathbb{E} Z_{t}^{2}]^{1/2}\bigg(e^{4(\alpha-\beta)t}\alpha^{4}\mathbb{E}\left(\int_{t-\delta}^{t}e^{-(\alpha-\beta)s}dM_{s}\right)^{4}
\nonumber
\\
&\qquad\qquad\qquad
+(e^{(\alpha-\beta)\delta}-1)^{4}e^{4(\alpha-\beta)(t-\delta)}\alpha^{4}
\mathbb{E}\left(\int_{0}^{t-\delta}e^{-(\alpha-\beta)s}dM_{s}\right)^{4}\bigg)^{1/2}.
\nonumber
\end{align*}
Now Proposition~\ref{prop:moments-noncritial} implies that $[\mathbb{E} Z_{t}^{2}]^{1/2} \le C n$. In conjunction with Proposition \ref{prop:ineq} and \eqref{eq:6},
we deduce from the above inequality that for $n$ sufficiently large,
\begin{align}
&\mathbb{E}\left[Z_{t}\left(Z_{t}-ne^{(\alpha-\beta)t}-Z_{t-\delta}
+ne^{(\alpha-\beta)(t-\delta)}\right)^{2}\right] \nonumber \\
&\leq Cn \cdot \left(C n^2 \delta^2 + C\delta^4 n^2 \right)^{\frac{1}{2}}\leq C n^{2}\delta.
\end{align}
Hence it follows that the second term on the right hand side of \eqref{MidStep} is also upper bounded by $C \delta^2$ for some constant $C.$
Therefore, we have proved (c) for the case $\alpha \ne \beta$.

We next proceed to prove (c) for the case $\alpha = \beta$. The proof is similar as for the case $\alpha \ne \beta$, so we only outline the key steps.
Given $Z_0=n$, we have $\tilde{Z}_{t}=\frac{Z_{t}-n}{\sqrt{n}}$. Then we have
\begin{align}\label{eq:2}
&\mathbb{E}\left[(\tilde{Z}_{t+\delta}-\tilde{Z}_{t})^{2}(\tilde{Z}_{t}-\tilde{Z}_{t-\delta})^{2}\right]
\\
&=\frac{1}{n^{2}}\mathbb{E}\left[(Z_{t+\delta}-Z_{t})^{2}(Z_{t}-Z_{t-\delta})^{2}\right]
\nonumber
\\
&=\frac{1}{n^{2}}\mathbb{E}\left[ (Z_{t}-Z_{t-\delta})^{2} \cdot \mathbb{E}[(Z_{t+\delta}-Z_{t})^{2} |Z_t] \right]
\nonumber
\\
&=\frac{\alpha^{2}\delta}{n^{2}}\mathbb{E}\left[Z_{t}(Z_{t}-Z_{t-\delta})^{2}\right],
\nonumber
\end{align}
where the second equality follows from the fact that $Z$ is Markovian, and the third equality follows from a similar argument as for \eqref{eq:1} in the case $\alpha \ne \beta$. To bound the last term in \eqref{eq:2} we
note that
\begin{align*}
\mathbb{E}\left[Z_{t}(Z_{t}-Z_{t-\delta})^{2}\right]
&=\mathbb{E}[Z_{t}^{3}]-2\mathbb{E}[Z_{t}^{2}Z_{t-\delta}]+\mathbb{E}[Z_{t}Z_{t-\delta}^{2}]
\\
&=\mathbb{E}[Z_{t-\delta}^{3}]+3\alpha^{2}\delta\mathbb{E}[Z_{t-\delta}^{2}]
+\frac{3}{2}\alpha^{4}\delta^{2}\mathbb{E}[Z_{t-\delta}]+\alpha^{3}\delta\mathbb{E}[Z_{t-\delta}]
\nonumber
\\
&\qquad
-2\mathbb{E}[Z_{t-\delta}^{3}]-2\alpha^{2}\delta\mathbb{E}[Z_{t-\delta}^{2}]
+\mathbb{E}[Z_{t-\delta}^{3}]
\nonumber\\
&=\alpha^{2}\delta\mathbb{E}[Z_{t-\delta}^{2}]
+\left(\frac{3}{2}\alpha^{4}\delta^{2}+\alpha^{3}\delta\right) \mathbb{E}[Z_{t-\delta}]
\nonumber \\
&=\alpha^{2}\delta(n^{2}+\alpha^{2}n(t-\delta))
+\left(\frac{3}{2}\alpha^{4}\delta^{2}+\alpha^{3}\delta\right)n,
\nonumber
\end{align*}
where the first equality follows from Proposition~\ref{prop:moments-noncritial} and tower property of conditional expectation, and the
last equality follows from Proposition~\ref{prop:moments-noncritial} and the fact that $Z_0=n$.
Hence,
\begin{equation*}
\mathbb{E}\left[(\tilde{Z}_{t+\delta}-\tilde{Z}_{t})^{2}(\tilde{Z}_{t}-\tilde{Z}_{t-\delta})^{2}\right]
=\alpha^{4}\delta^{2}\left[1+\frac{\alpha^{2}}{n}(t-\delta)+\frac{3}{2n}\alpha^{2}\delta+\frac{\alpha}{n}\right]
\leq C \delta^{2},
\end{equation*}
where $C$ is a positive constant that is independent of $n$.
Hence we have also established (c) when $\alpha = \beta$. The proof is therefore complete.

\subsubsection{FCLT for $N$ when $\mu=0$}
In this section we prove the weak convergence of the sequence of re--normalized Hawkes processes \eqref{eq:FCLT-N} when $\mu=0$.
Let us recall that the intensity process $Z_{t}$ satisfies the dynamics
\begin{equation*}
dZ_{t}=-\beta Z_{t}dt+\alpha dN_{t},\qquad Z_{0}=n,
\end{equation*}
where $N_{t}$ is a simple point process with intensity $Z_{t-}$. We can express
the jump process $N_{t}=N(0,t]$ in terms of the intensity process $Z_{t}$ in the following way,
\begin{equation*}
N_{t}=\frac{Z_{t}-Z_{0}}{\alpha}+\frac{\beta}{\alpha}\int_{0}^{t}Z_{s}ds.
\end{equation*}
This immediately yields that
\begin{align}\label{eq:N-convg}
\frac{N_{t}-n\int_{0}^{t}e^{(\alpha-\beta)s}ds}{\sqrt{n}}
&=\frac{1}{\sqrt{n}}\left[\frac{Z_{t}-n}{\alpha}+\frac{\beta}{\alpha}\int_{0}^{t}Z_{s}ds
-n\int_{0}^{t}e^{(\alpha-\beta)s}ds\right]
\nonumber\\
&=\frac{1}{\sqrt{n}}\left[\frac{Z_{t}-ne^{(\alpha-\beta)t}}{\alpha}
+\frac{\beta}{\alpha}\int_{0}^{t}(Z_{s}-ne^{(\alpha-\beta)s})ds\right]
\end{align}
Note that we have shown that for any $T>0$, as $n\rightarrow\infty$,
\begin{align*}
\left\{\frac{Z_{t}-ne^{(\alpha-\beta)t}}{\sqrt{n}}: t \in [0,T] \right\}\rightarrow G,
\end{align*}
weakly on $D[0,T]$ equipped with Skorohod $J_1$ topology,
where $G$ is a centered Gaussian process. Thus, using  \eqref{eq:N-convg}, we conclude from \cite[Theorem~2.2]{kurtz1991weak} that as $n\rightarrow\infty$,
\begin{align}
\left\{ \frac{N_{t}-n\int_{0}^{t}e^{(\alpha-\beta)s}ds}{\sqrt{n}} : t \in [0,T] \right\} \rightarrow H,
\nonumber
\end{align}
weakly on $D[0,T]$, where $H$ is given in \eqref{eq:H}.

\subsubsection{FCLT for $Z$ and $N$ when $\mu>0$}

In this section we prove Theorem~\ref{thm:CLT} for the case $\mu>0$ using the observation in Section~\ref{sec:mu-relation}.

Decompose $N_{t}=N_{t}^{(0)}+N_{t}^{(1)}$.  Note for any $T>0$, we have $\sup_{0\leq t\leq T}N_{t}^{(1)}$ is finite
almost surely and independent
of the parameter $n$. This implies that as $n\rightarrow\infty$
\begin{equation} \label{eq:N1}
\sup_{0\leq t\leq T} \frac{N_{t}^{(1)}}{\sqrt{n}} \rightarrow 0, \quad \text{almost surely}.
\end{equation}
Since we have established FCLT for $N_{t}^{(0)}$, the result for $N$ then readily follows from \eqref{eq:N1}.

Similarly, we can decompose $Z_{t}=Z_{t}^{(0)}+Z_{t}^{(1)}$ and note
that $Z_{t}^{(1)}$ is independent of $n$ and hence $\sup_{0\leq t\leq T}Z_{t}^{(1)}/\sqrt{n} \le \alpha \sup_{0\leq t\leq T}N_{t}^{(1)}/\sqrt{n}\rightarrow 0$ almost surely as $n\rightarrow\infty$.
Hence the FCLT for $Z$ follows.

\subsection{Proof of Theorem 3}
We prove Theorem~\ref{thm:rescale-LLN} in this section. We also first prove Theorem~\ref{thm:rescale-LLN} when $\mu=0$, and then prove it when $\mu>0$ using the observation in Section~\ref{sec:mu-relation}.

\subsubsection{Proof of Theorem~3 when $\mu=0$}

\begin{proof} [Proof of part (i): super--critical case]

We first show \eqref{eq:Zlargetime}.
Observe first from \eqref{MartingaleRep} that for any $n\in\mathbb{N}$,
\begin{equation*}
\frac{Z_{s \tau_{n}}-n^{1+s}}{n^{1+s}}=\frac{\alpha}{n}\int_{0}^{s \tau_{n}}e^{-(\alpha-\beta)u}dM_{u}
\end{equation*}
is a martingale. Therefore, Doob's martingale inequality implies that for any $\epsilon>0$,
\begin{eqnarray*}
\mathbb{P}\left(\sup_{0\leq s\leq T}\left|\frac{Z_{s \tau_{n}}-n^{1+s}}{n^{1+s}}\right|\geq\epsilon\right)
&=&\mathbb{P}\left(\sup_{0\leq s\leq T}\left|\int_{0}^{s \tau_{n}}e^{-(\alpha-\beta)u}dM_{u}\right|
\geq\frac{n\epsilon}{\alpha}\right)
\nonumber
\\
&\leq&\frac{C\alpha^{4}}{n^{4}\epsilon^{4}}
\mathbb{E}\left[\left(\int_{0}^{T \tau_{n}}e^{-(\alpha-\beta)u}dM_{u}\right)^{4}\right].
\nonumber
\end{eqnarray*}
Now inequality \eqref{ineq:mart} implies that for $\alpha> \beta$ and $\tau_{n}=\frac{\log n}{\alpha-\beta}$,
\begin{align}
\mathbb{E}\left[\left(\int_{0}^{T \tau_{n}}e^{-(\alpha-\beta)u}dM_{u}\right)^{4}\right]
\le C T \tau_n \cdot \left( n^2 (1-n^{-2T}) + n (1-n^{-2T}) + n (1-n^{-3T}) \right). \nonumber
\end{align}
Thus we obtain for $n$ sufficiently large,
\begin{align}
\mathbb{P}\left(\sup_{0\leq s\leq T}\left|\frac{Z_{s \tau_{n}}-n^{1+s}}{n^{1+s}}\right|\geq\epsilon\right)
\leq\frac{C T \log n}{n^{2}\epsilon^{4}}.\nonumber
\end{align}
Hence, by Borel-Cantelli lemma,
\begin{equation*}
\sup_{0\leq s\leq T}\left|\frac{Z_{s \tau_{n}}}{n^{1+s}}-1\right|\rightarrow 0,
\end{equation*}
almost surely as $n\rightarrow\infty$.

We next prove \eqref{eq:Nlargetime}. Recall that the intensity process $Z_{t}$ satisfies the dynamics
\begin{equation*}
dZ_{t}=-\beta Z_{t}dt+\alpha dN_{t},\qquad Z_{0}=n,
\end{equation*}
where $N_{t}$ is a simple point process with intensity $Z_{t-}$ (since $\mu=0$). We can therefore express
the jump process $N_{t}=N(0,t]$ in terms of the intensity process $Z_{t}$ in the following way,
\begin{equation*}
N_{t}=\frac{Z_{t}-Z_{0}}{\alpha}+\frac{\beta}{\alpha}\int_{0}^{t}Z_{s}ds.
\end{equation*}
This implies that for any $t>0$,
\begin{equation*} \label{eq:N-Z-centered}
N_{t}-n\int_{0}^{t}e^{(\alpha-\beta)s}ds
=\frac{Z_{t}-ne^{(\alpha-\beta)t}}{\alpha}
+\frac{\beta}{\alpha}\int_{0}^{t}(Z_{s}-ne^{(\alpha-\beta)s})ds.
\end{equation*}
Hence, for $\alpha> \beta$ and $\tau_n = \frac{\log n}{\alpha- \beta}$, we obtain
\begin{align*}
\left|\frac{N_{s\tau_{n}}-\frac{n^{1+s}-n}{\alpha-\beta}}{n^{1+s}}\right|
&=\left|\frac{Z_{s \tau_{n}}-n^{1+s}}{\alpha n^{1+s}}
+\frac{\beta}{\alpha} \tau_{n}\int_{0}^{s}n^{u-s}\frac{Z_{u \tau_{n}}-n^{1+u}}{n^{1+u}}du\right|
\\
&\leq\left|\frac{Z_{s\tau_{n}}-n^{1+s}}{\alpha n^{1+s}}\right|+\frac{\beta}{\alpha} \tau_n
\sup_{0\leq u\leq s}\left|\frac{Z_{u\tau_{n}}-n^{1+u}}{n^{1+u}}\right| \cdot \int_{0}^{s}n^{u-s} du .
\nonumber
\\
&\leq\frac{1}{\alpha} \cdot\left|\frac{Z_{s\tau_{n}}-n^{1+s}}{ n^{1+s}}\right|+\frac{\beta}{\alpha (\alpha- \beta)}
\sup_{0\leq u\leq s}\left|\frac{Z_{u\tau_{n}}-n^{1+u}}{n^{1+u}}\right|.
\nonumber
\end{align*}
Then \eqref{eq:Zlargetime} implies the desired result \eqref{eq:Nlargetime}.
\end{proof}

\begin{proof} [Proof of part (ii): sub--critical case]

The proof is similar as in the super--critical case, so we only outline the key steps.

We first show \eqref{eq:Zlargetime-sub}.
For $\alpha< \beta$, we observe that for $0\leq s\leq T<1$,
\begin{equation*}
\frac{Z_{st_{n}}-n^{1-s}}{n^{1-s}}=\frac{\alpha}{n}\int_{0}^{st_{n}}e^{-(\alpha-\beta)s}dM_{s}
\end{equation*}
is a martingale. Thus we have, for $t_{n}=\frac{\log n}{\beta -\alpha}$ and any $\epsilon>0$,
\begin{eqnarray*}
\mathbb{P}\left(\sup_{0\leq s\leq T}\left|\frac{Z_{s t_{n}}-n^{1-s}}{n^{1-s}}\right|\geq\epsilon\right)
&\leq& \frac{C}{n^{4}\epsilon^{4}}
\mathbb{E}\left[\left(\int_{0}^{T t_{n}}e^{-(\alpha-\beta)u}dM_{u}\right)^{4}\right]
\nonumber\\
&\leq & \frac{C}{n^{4}\epsilon^{4}} T t_n \cdot \left( n^2 (n^{2T}-1) + n (n^{2T}-1) + n (n^{3T}-1) \right),
\nonumber
\end{eqnarray*}
where the last inequality follows from \eqref{ineq:mart}. Therefore, for $t_{n}=\frac{\log n}{\beta -\alpha}$, we obtain for $n$ sufficiently large,
\begin{align*}
\mathbb{P}\left(\sup_{0\leq s\leq T}\left|\frac{Z_{st_{n}}-n^{1-s}}{n^{1-s}}\right|\geq\epsilon\right)
\leq
\frac{C T {\log n} } {n^{4}\epsilon^{4}}
\left[{n^{2}n^{2T}}+{n\cdot n^{3T}}\right],
\nonumber
\end{align*}
where the right-hand side of the above inequality goes to $0$ as $n\rightarrow\infty$ for any $T<1$. Thus we have $\sup_{0\leq s\leq T}\left|\frac{Z_{s t_{n}}-n^{1-s}}{n^{1-s}}\right|$
converges to zero in probability for any $T<1$.
Furthermore, if $T<\frac{1}{2}$, we have \[\sum_{n=1}^{\infty} \frac{{\log n} } {n^{4}\epsilon^{4}}
\left[{n^{2}n^{2T}}+{n\cdot n^{3T}}\right] < +\infty.\] Hence the almost sure convergence follows from Borel-Cantelli lemma.

We next prove \eqref{eq:Nlargetime-sub}. When $\alpha < \beta$ and $t_n = \frac{\log n}{\beta - \alpha}$, we have
\begin{align}
\left|\frac{N_{st_{n}}-\frac{1}{\beta-\alpha}(n-n^{1-s})}{n}\right|
&=\left|\frac{Z_{st_{n}}-n^{1-s}}{\alpha n}
+\frac{\beta}{\alpha}t_{n}\int_{0}^{s}n^{-u}\frac{Z_{ut_{n}}-n^{1-u}}{n^{1-u}}du\right| \nonumber
\\
&\leq n^{-s}\left|\frac{Z_{st_{n}}-n^{1-s}}{\alpha n^{1-s}}\right|+\frac{\beta}{\alpha (\beta - \alpha)}\sup_{0\leq u\leq s}\left|\frac{Z_{ut_{n}}-n^{1-u}}{n^{1-u}}\right|.
\nonumber
\end{align}
The desired result \eqref{eq:Nlargetime-sub} then follows.
\end{proof}

\subsubsection{Proof of Theorem~3 when $\mu>0$}

In this section we prove Theorem~3 for the case $\mu>0$ using the observation in Section~\ref{sec:mu-relation}.

\textit{Proof of \eqref{eq:Zlargetime}}:
Note that $Z_{t}^{(1)}\leq\alpha N_{t}^{(1)}$. 
It was in Section 5.4. in Zhu \cite{ZhuThesis} 
that $\lim_{t\rightarrow\infty}\frac{1}{t}\log N_{t}^{(1)}=\alpha-\beta$ almost surely
for $\alpha>\beta>0$. Thus, for any $0<T<1$,
\begin{equation*}
\sup_{0\leq s\leq T}\frac{Z_{s\tau_{n}}^{(1)}}{n^{1+s}}
\leq\alpha\sup_{0\leq s\leq T}\frac{N_{s\tau_{n}}^{(1)}}{n^{1+s}}\leq\frac{N_{T\tau_{n}}^{(1)}}{n}\rightarrow 0,
\end{equation*}
almost surely as $n\rightarrow\infty$ since $\tau_{n}=\frac{\log n}{\alpha- \beta}$. Then \eqref{eq:Zlargetime} follows from the decomposition
$Z_t= Z_{t}^{(0)} +Z_{t}^{(1)}$ and the result proved in the previous section.

\textit{Proof of \eqref{eq:Nlargetime}}: Similarly, it follows from the decomposition $N_t= N_{t}^{(0)} +N_{t}^{(1)}$
and that for any $0<T<1$,
\begin{equation*}
\sup_{0\leq s\leq T}\frac{N_{s \tau_{n}}^{(1)}}{n^{1+s}}
\leq\frac{N_{T \tau_{n}}^{(1)}}{n}\rightarrow 0,
\end{equation*}
almost surely as $n\rightarrow\infty$ since $ \tau_{n}=\frac{\log n}{\alpha-\beta}$.

\textit{Proof of \eqref{eq:Zlargetime-sub}}:
Note that for $\beta>\alpha>0$, $\frac{N_{t}^{(1)}}{t}\rightarrow\frac{\mu}{1-\frac{\alpha}{\beta}}$ almost surely
as $t\rightarrow\infty$. Therefore, for $t_{n}=\frac{\log n}{\beta-\alpha}$ we deduce that for any $0<T<1$,
\begin{equation*}
\sup_{0\leq s\leq T}\frac{Z_{st_{n}}^{(1)}}{n^{1-s}}
\leq\frac{\alpha N_{Tt_{n}}^{(1)}}{n^{1-T}}\rightarrow 0,
\end{equation*}
almost surely as $n\rightarrow\infty$.

\textit{Proof of \eqref{eq:Nlargetime-sub}}: Similarly, for $t_{n}=\frac{\log n}{\beta-\alpha}$, we have
\begin{equation*}
\sup_{0\leq s\leq T}\frac{N_{st_{n}}^{(1)}}{n}=\frac{N_{Tt_{n}}^{(1)}}{n}\rightarrow 0,
\end{equation*}
almost surely as $n\rightarrow\infty$ since  $\frac{N_{t}^{(1)}}{t}\rightarrow\frac{\mu}{1-\frac{\alpha}{\beta}}$
almost surely as $t\rightarrow\infty$.
\end{proof}

\subsection{Proof of Theorem~4}

\subsubsection{Proof of critical and nearly--critical cases}

\begin{proof}
The proof is based on diffusion approximations. In particular, we apply Theorem~4.1 in Ethier and Kurtz \cite[Chapter~7]{EthierKurtz} and verify their conditions~(4.1)--(4.7). Fix $Z_0=n$. Define for $t \ge 0$ and $\mu \ge 0$,
\begin{eqnarray} \label{eq:XBA}
\mathbb{X}^n_t := \frac{Z_{nt}}{n}, \quad \mathbb{B}^n_t := \int_{0}^{t} (\alpha_n \mu + \gamma \mathbb{X}^n_s) ds, \quad
\mathbb{A}^n_t := \alpha_n^2 \int_{0}^{t} \left(\mathbb{X}^n_s + \frac{\mu}{n}\right) ds.
\end{eqnarray}
One readily checks that the two processes
\begin{eqnarray*}
\left\{\mathbb{X}^n_t -\mathbb{B}^n_t: t \in[0, T] \right\},
\quad \text{and} \quad \left\{(\mathbb{X}^n_t -\mathbb{B}^n_t)^2- \mathbb{A}^n_t: t \in[0, T]\right\}
\end{eqnarray*}
are martingales with respect to the filtration generated by $\mathbb{X}^n$. Since $Z$ can only make jumps of size $\alpha_n$, we deduce that $\mathbb{X}^n$ can only
make jumps of size $\frac{\alpha_n}{n}$. This implies condition~(4.3) of Theorem~4.1 in \cite[Chapter~7]{EthierKurtz} holds. Conditions (4.4) and (4.5) hold trivially since
$\mathbb{B}^n$ and $\mathbb{A}^n$ have continuous sample paths. In view of \eqref{eq:XBA} and $\lim_{n \rightarrow \infty}\alpha_n = \beta$, Conditions (4.6) and (4.7) also hold if we set $b(x) = \beta \mu +\gamma x$ and $a(x)= \beta^2 x$. Now if we define
\[dX_{t}=(\beta \mu + \gamma X_{t})dt+\beta\sqrt{X_{t}}dB_{t},\qquad X_{0}=1\]
where $B$ is a standard Brownian motion, then this stochastic differential equation (SDE) has a pathwise unique strong solution. Since $\mathbb{X}^n_0 \equiv 1$
for all $n$, we deduce that $\mathbb{X}^n$ converges weakly to $X$ by Theorem~4.1 in \cite[Chapter~7]{EthierKurtz} and the equivalence of SDEs and martingale problems (see, e.g., \cite{kurtz2011} for details).

We next establish the weak convergence of the sequence of re--normalized Hawkes processes \[\left\{\frac{N_{tn}}{n^2}: t \in [0,T] \right\}.\]
Since for $\mu \ge 0$, the point process $N$ has intensity $\mu+Z_{t-}$ at time $t$, where
$dZ_{t}=-\beta Z_{t}dt+\alpha_{n} dN_{t}$ for $Z_{0} = n$. So we still have
\begin{equation*}
N_{t}=\frac{Z_{t}-Z_{0}}{\alpha_n}+\frac{\beta}{\alpha_n}\int_{0}^{t}Z_{s}ds,
\end{equation*}
which yields
\begin{equation*}
\frac{N_{tn}}{n^2} =\frac{Z_{tn}-Z_{0}}{n^2 \cdot \alpha_n}+\frac{\beta}{\alpha_n} \int_{0}^{t}\frac{Z_{sn}}{n} ds.
\end{equation*}
Since $\lim_{n \rightarrow \infty} {\alpha_n} = \beta$ and we have the weak convergence of the sequence of processes $\left\{\frac{Z_{tn}}{n}: t \in [0,T] \right\}$, the result then follows from \cite[Theorem~2.2]{kurtz1991weak}.
\end{proof}

\subsubsection{Proof of super--critical case when $\mu=0$}

\begin{proof} The proof is based on Aldous's result on weak convergence of a sequence of martingales to a continuous
martingale limit \cite{Aldous1989}.

We first note that in the super-critical case $\alpha>\beta>0$ where $\tau_{n}:=\frac{\log n}{\alpha-\beta}$, we have
\begin{equation*}
Z_{s\tau_n}-n^{1+s}
=\alpha n^{s}\int_{0}^{s\tau_n}e^{-(\alpha-\beta)u}dM_{u}.
\end{equation*}
Therefore, $\left\{\frac{Z_{s\tau_n}-n^{1+s}}{n^{\frac{1}{2}+s}}: s \in [0, T]\right\}$ is a martingale
for any fixed $n$. Moreover, recall that $\langle M\rangle_{t}=\int_{0}^{t}Z_{s}ds$
and it is straightforward to compute that
\begin{align*}
\sup_{n\in\mathbb{N}}\mathbb{E}\left[\left(
\frac{Z_{s\tau_n}-n^{1+s}}{n^{\frac{1}{2}+s}}\right)^{2}
\right]
&=\sup_{n\in\mathbb{N}}\frac{\alpha^{2}}{n}\mathbb{E}\left[\int_{0}^{s\tau_n}
e^{-2(\alpha-\beta)u}Z_{u}du\right]
\\
&=\sup_{n\in\mathbb{N}}\frac{\alpha^{2}}{n}\int_{0}^{s\tau_n}e^{-(\alpha-\beta)u}Z_{0}du
\nonumber
\\
&=\sup_{n\in\mathbb{N}}\frac{\alpha^{2}}{\alpha-\beta}\left(1-\frac{1}{n^{s}}\right)
\nonumber
\\
&=\frac{\alpha^{2}}{\alpha-\beta}.
\nonumber
\end{align*}
Therefore, for each $s$, $\left\{\frac{Z_{s\tau_n}-n^{1+s}}{n^{\frac{1}{2}+s}}: n \ge 1 \right\}$
is uniformly integrable.

We next establish the convergence of finite-dimensional distributions. To this end, we first show for fixed $s>0$, the sequence of moment generating functions
\begin{equation*}
\mathbb{E}\left[e^{-\theta\frac{Z_{s \tau_n}-n^{1+s} }{\sqrt{n^{1+2s}}}}\bigg|Z_{0}=n\right]
\end{equation*}
converges when $n\rightarrow \infty$. 
First, we notice from \eqref{eq:u} that for any $\theta\in\mathbb{R}$,
when $n$ is sufficiently large, we have
\begin{equation} \label{eq:mgf1-rescaled}
\mathbb{E}\left[e^{-\theta\frac{Z_{s\tau_n}-n^{1+s} }{\sqrt{n^{1+2s}}}}\bigg|Z_{0}=n\right]
=\exp \left({ A\left(s\tau_n, -\frac{\theta}{\sqrt{n^{1+2s}}} \right) \cdot n+ \sqrt{n}\theta } \right),
\end{equation}
where $A$ solves the ODE in \eqref{eq:A}. To see this, it suffices to show that for any $\theta\in\mathbb{R}$, $\frac{|\theta|}{\sqrt{n^{1+2s}}}<\theta_{c}(s\tau_{n})$
for sufficiently large $n$,
where $\theta_{c}(\cdot)$ is defined in \eqref{thetacEqn}. 
Note that in the super--critical case $\alpha>\beta$, by the definition of $\theta_{c}(\cdot)$
in \eqref{thetacEqn}, we can show
that for any $\epsilon>0$, 
we have $\theta_{c}(t)\geq e^{-(\alpha-\beta+\epsilon)t}$ for sufficiently large $t>0$.
To see this, we recall the definition of $\theta_{c}(\cdot)$ in \eqref{thetacEqn}
and in the super--critical case $\alpha>\beta$, $-\beta A+e^{\alpha A}-1\geq 0$
for every $A\geq 0$ and it is zero if and only if $A=0$. Thus, $\theta_{c}(t)\rightarrow 0$
as $t\rightarrow\infty$. For any $\epsilon>0$, there exists $\eta>0$ sufficiently small, so that
for any $0\leq A\leq\eta$, $-\beta A+e^{\alpha A}-1\leq\left(\alpha-\beta+\frac{1}{2}\epsilon\right)A$.
Therefore, by the definition of $\theta_{c}(t)$, we have for sufficiently large $t$, $\theta_{c}(t)<\eta$ and
\begin{align*}
t&=\int_{\theta_{c}(t)}^{\eta}\frac{dA}{-\beta A+e^{\alpha A}-1}
+\int_{\eta}^{\infty}\frac{dA}{-\beta A+e^{\alpha A}-1}
\\
&\geq
\int_{\theta_{c}(t)}^{\eta}\frac{dA}{\left(\alpha-\beta+\frac{1}{2}\epsilon\right)A}
+\int_{\eta}^{\infty}\frac{dA}{-\beta A+e^{\alpha A}-1}
\\
&=-\frac{1}{\alpha-\beta+\frac{1}{2}\epsilon}\log(\theta_{c}(t))+\frac{1}{\alpha-\beta+\frac{1}{2}\epsilon}\log(\eta)
+\int_{\eta}^{\infty}\frac{dA}{-\beta A+e^{\alpha A}-1},
\end{align*}
which implies that
\begin{equation}
\theta_{c}(t)\geq
\eta e^{(\alpha-\beta+\frac{1}{2}\epsilon)\int_{\eta}^{\infty}\frac{dA}{-\beta A+e^{\alpha A}-1}}
e^{-(\alpha-\beta+\frac{1}{2}\epsilon)t}
\geq e^{-(\alpha-\beta+\epsilon)t},
\end{equation}
for any sufficiently large $t$. 
Thus, $\theta_{c}(s\tau_{n})\geq n^{-\frac{\alpha-\beta+\epsilon}{\alpha-\beta}s}>\frac{|\theta|}{\sqrt{n^{1+2s}}}$
for any sufficiently large $n$ and small $\epsilon$. 
Second, we will show below that
\begin{eqnarray} \label{eq:A-asy}
\lim_{n \rightarrow \infty}  \left({ A\left(s\tau_n, -\frac{\theta}{\sqrt{n^{1+2s}}} \right) \cdot n+ \sqrt{n}\theta } \right) =  \frac{1}{2} \frac{\alpha^{2}}{\alpha-\beta} {\theta^2} ,
\end{eqnarray}
and this implies that for fixed $s >0$,
\begin{equation*}
\frac{Z_{s\tau_n}-n^{1+s}}{\sqrt{n^{1+2s}}} \rightarrow \xi, \quad \text{weakly as $n \rightarrow \infty$},
\end{equation*}
where $\xi$ is a normal random variable with mean zero and variance $\frac{\alpha^{2}}{\alpha-\beta}$.

To establish \eqref{eq:A-asy}, we rely on Gronwall's inequality to obtain estimates of $A$.
Write $g(x) = e^{\alpha x} - \alpha x -1$.
Then given any small $\epsilon, \eta>0$ with $\epsilon < \frac{\alpha^2}{2}$, there exists some $\delta>0$ such that $ \left( \frac{1}{2} \alpha^2 - \epsilon \right) x^2 \le g(x) \le \left( \frac{1}{2} \alpha^2 + \epsilon \right) x^2$ and $g(x) \le  \eta |x|$ when $|x| \le \delta$. Recall from \eqref{eq:A} that $A$ solves the ODE:
\begin{align}
A'\left(t\right) =(\alpha - \beta)A\left(t\right)+ g\left(A\left(t\right)\right). \label{eq:ODE-A1}
\end{align}
Suppose that for $n$ large, we have $\left|A\left(t; -\frac{\theta}{\sqrt{n^{1+2s}}} \right)\right| \le \delta$ for all $t \in [0, \tau_n s]$. Then we deduce that for $t \in [0, \tau_n s]$,
\begin{eqnarray} \label{eq:A-ineq}
A'\left(t; -\frac{\theta}{\sqrt{n^{1+2s}}} \right) \le (\alpha - \beta) \cdot A\left(t;-\frac{\theta}{\sqrt{n^{1+2s}}} \right) +  \left( \frac{1}{2} \alpha^2 + \epsilon \right) A^2\left(t;-\frac{\theta}{\sqrt{n^{1+2s}}} \right).
\end{eqnarray}
The solution to the Bernoulli equation
\begin{eqnarray} \label{eq:y-eq}
y'(t) &=& (\alpha - \beta)y(t) +  \left( \frac{1}{2} \alpha^2 + \epsilon \right)  y^2(t),
\end{eqnarray}
is given by
\[y(t) = \left( \left(\frac{1}{y(0)} + \frac{\frac{1}{2} \alpha^2 + \epsilon }{\alpha - \beta}\right) \cdot e^{(\beta - \alpha) t} - \frac{\frac{1}{2} \alpha^2 + \epsilon}{\alpha - \beta} \right)^{-1}.\]
It readily follows that given $y(0)=-\frac{\theta}{\sqrt{n^{1+2s}}}$, we have for $0<s<\frac{1}{2}$,
\begin{eqnarray} \label{eq:y-asy}
\lim_{n \rightarrow \infty}  \left({ y\left(s\tau_n \right) \cdot n+ \sqrt{n}\theta } \right) =   \frac{\frac{1}{2} \alpha^2 + \epsilon}{\alpha-\beta} {\theta^2} ,
\end{eqnarray}
In addition, we have $|y(t)| \le \delta$ for all $t \in [0, \tau_n s]$ when $n$ is large. Note that the quadratic function $(\alpha - \beta)y+  \left( \frac{1}{2} \alpha^2 + \epsilon \right)  y^2$ is Lipschitz continuous in $y$ when $|y| \le \delta$. Then using Gronwall's inequality for nonlinear ODEs, we can infer from \eqref{eq:A-ineq} and \eqref{eq:y-eq} that
\[A\left(t;-\frac{\theta}{\sqrt{n^{1+2s}}} \right) \le y(t), \quad \text{for all $t \in [0, \tau_n s]$.} \]
On combining \eqref{eq:y-asy}, we then deduce that
\begin{eqnarray}
\limsup_{n \rightarrow \infty}  \left({ A\left(s\tau_n, -\frac{\theta}{\sqrt{n^{1+2s}}} \right) \cdot n+ \sqrt{n}\theta } \right) \le \frac{\frac{1}{2} \alpha^2 + \epsilon}{\alpha-\beta} {\theta^2}.
\end{eqnarray}
A similar argument yields the lower bound:
\begin{eqnarray}
\liminf_{n \rightarrow \infty}  \left({ A\left(s\tau_n, -\frac{\theta}{\sqrt{n^{1+2s}}} \right) \cdot n+ \sqrt{n}\theta } \right) \ge \frac{\frac{1}{2} \alpha^2 - \epsilon}{\alpha-\beta} {\theta^2}.
\end{eqnarray}
Letting $\epsilon \rightarrow 0$, we obtain \eqref{eq:A-asy}.

It only remains to show that for $n$ large, we have $\left|A\left(t; -\frac{\theta}{\sqrt{n^{1+2s}}} \right)\right| \le \delta$ for all $t \in [0, \tau_n s]$.
To this end, we first define for fixed $s$ and $n$,
\[c_{n}(s):= \sup \left\{u \ge0 : \left|A \left(t;  -\frac{\theta}{\sqrt{n^{1+2s}}} \right)\right| \le \delta, \quad \text{for all $ t \le u$} \right\}.\]
Note that $\left|A\left(0; -\frac{\theta}{\sqrt{n^{1+2s}}} \right)\right| =\left| -\frac{\theta}{\sqrt{n^{1+2s}}} \right| < \delta $ for all large $n$. Thus $c_{n}(s)>0$ for all large $n$.
It suffices to show $\tau_n s \le c_{n}(s)$ for all sufficiently large $n$. If $c_{n}(s)=\infty$, then the inequality $\tau_{n}s\leq c_{n}(s)$ holds automatically.
If $c_{n}(s)<\infty$, we proceed as follows.
We note from the ODE in \eqref{eq:ODE-A1} that
\[A\left(t; -\frac{\theta}{\sqrt{n^{1+2s}}} \right) = A\left(0;-\frac{\theta}{\sqrt{n^{1+2s}}} \right) \cdot e^{(\alpha - \beta)t}
+ \int_{0}^t e^{(\alpha - \beta) (t-u)} g\left(A\left(u; -\frac{\theta}{\sqrt{n^{1+2s}}} \right)\right) du. \]
Note that $\left|A\left(t; -\frac{\theta}{\sqrt{n^{1+2s}}} \right)\right| \le \delta$ for all $ t \in [0, c_{n}(s)]$. In addition, $g(x) \le  \eta |x|$ when $|x| \le \delta$.
Hence we obtain for $t \in [0, c_{n}(s)]$,
\[ e^{-(\alpha - \beta)t} \left|A\left(t; -\frac{\theta}{\sqrt{n^{1+2s}}} \right)\right| \le \left|A\left(0; -\frac{\theta}{\sqrt{n^{1+2s}}} \right)\right|
+ \eta \int_{0}^t e^{-(\alpha - \beta)u} \left|A\left(u; -\frac{\theta}{\sqrt{n^{1+2s}}} \right)\right| du.  \]
Gronwall's inequality then implies
\begin{eqnarray} \label{eq:bound-A:2}
\left|A\left(t; -\frac{\theta}{\sqrt{n^{1+2s}}} \right)\right| \le \frac{|\theta|}{\sqrt{n^{1+2s}}} \cdot e^{(\alpha - \beta)t} \cdot e^{\eta t}, \quad \text{for $t \in [0,c_{n}(s)]$}.
\end{eqnarray}
For $c_{n}(s)<\infty$, letting $t=c_{n}(s)$ in \eqref{eq:bound-A:2}, we obtain from the definition of $c_{n}(s)$ that
\begin{equation}
\delta=\left|A\left(c_{n}(s); -\frac{\theta}{\sqrt{n^{1+2s}}} \right)\right|
\leq
\frac{|\theta|}{\sqrt{n^{1+2s}}} \cdot e^{(\alpha - \beta)c_{n}(s)} \cdot e^{\eta c_{n}(s)},
\end{equation}
which implies that
\begin{equation}
c_{n}(s)\geq\frac{\log(\delta/|\theta|)}{\alpha-\beta+\eta}
+\frac{\frac{1}{2}+s}{\alpha-\beta+\eta}\log n.
\end{equation}
Now given any $\eta<\frac{\alpha-\beta}{2T}$, and any $s \in [0, T]$, there exists a positive integer $K$ that is independent of $s$ such that for all $n>K$, we have
\[\frac{\log(\delta/|\theta|)}{\alpha-\beta+\eta}
+\frac{\frac{1}{2}+s}{\alpha-\beta+\eta}\log n\geq \tau_{n}s=\frac{\log n}{\alpha-\beta}s.\]
Hence we have $\tau_{n}s\le c_{n}(s)$ for all sufficiently large $n$
and thus we have proved that $\left|A\left(t; -\frac{\theta}{\sqrt{n^{1+2s}}} \right)\right| \le \delta$ for all $t \in [0, \tau_n s]$.

After establishing the convergence of one-dimensional marginal distributions, we proceed to consider the dimension of two. The general case of the convergence of finite--dimensional distributions follows similarly and the proof is omitted. Fix $\theta_1, \theta_2 \in\mathbb{R}$.
For any $0<u<v \le T$, and for sufficiently large $n$, the moment generating function of the random vector
$\left(\frac{Z_{u\tau_n}-n^{1+u}}{\sqrt{n^{1+2u}}}, \frac{Z_{v\tau_n}-n^{1+v}}{\sqrt{n^{1+2v}}} \right)$ is
given by
\begin{eqnarray} \label{eq:mgf-2d-rescaled}
\mathbb{E}\left[e^{-\frac{\theta_{1}}{\sqrt{n}n^{u}}Z_{u\tau_n}
-\frac{\theta_{2}}{\sqrt{n}n^{v}}Z_{v\tau_n}}\bigg|Z_{0}=n\right] \cdot e^{\theta_{1}\sqrt{n}
+\theta_{2}\sqrt{n}}.
\end{eqnarray}
By Cauchy-Schwarz inequality,
\begin{align*}
&\mathbb{E}\left[e^{-\frac{\theta_{1}}{\sqrt{n}n^{u}}Z_{u\tau_n}
-\frac{\theta_{2}}{\sqrt{n}n^{v}}Z_{v\tau_n}}\bigg|Z_{0}=n\right]
\\
&\leq
\left(\mathbb{E}\left[e^{-\frac{2\theta_{1}}{\sqrt{n}n^{u}}Z_{u\tau_n}}\bigg|Z_{0}=n\right]\right)^{\frac{1}{2}}
\left(\mathbb{E}\left[e^{-\frac{2\theta_{2}}{\sqrt{n}n^{v}}Z_{v\tau_n}}\bigg|Z_{0}=n\right]\right)^{\frac{1}{2}}
<\infty,
\end{align*}
for any sufficiently large $n$. It can be directly computed that
\begin{eqnarray*}
\lefteqn{\mathbb{E}\left[e^{-\frac{\theta_{1}}{\sqrt{n}n^{u}}Z_{u\tau_n}
-\frac{\theta_{2}}{\sqrt{n}n^{v}}Z_{v\tau_n}}\bigg|Z_{0}=n\right]} \\
& =& \mathbb{E}\left[e^{-\frac{\theta_{1}}{\sqrt{n}n^{u}}Z_{u\tau_n}}
\cdot \mathbb{E} \left[e^{- \frac{\theta_{2}}{\sqrt{n}n^{v}}Z_{v\tau_n}} \bigg| Z_{u \tau_n} \right]\Big|Z_{0}=n \right] \\
& =& \mathbb{E}\left[e^{-\frac{\theta_{1}}{\sqrt{n}n^{u}}Z_{u\tau_n}}
\cdot \left[ e^{ A((v-u) \tau_n , - \frac{\theta_{2}}{\sqrt{n}n^{v}} ) \cdot Z_{u \tau_n} } \right]\Big|Z_{0}=n \right],
\end{eqnarray*}
hence the moment generating function in \eqref{eq:mgf-2d-rescaled} becomes,
\begin{equation*}
\mathbb{E}\left[e^{- \hat \theta_1^{(n)} \frac{Z_{u\tau_n}-n^{1+u} }{\sqrt{n^{1+2u}}} }\bigg|Z_{0}=n\right] \cdot \exp(\Gamma^n),
\end{equation*}
where
\begin{eqnarray*}
\hat \theta_1^{(n)} &:=&\theta_1 - \sqrt{n^{1+ 2u}}  \cdot A \left((v-u) \tau_n , - \frac{\theta_{2}}{\sqrt{n}n^{v}} \right),\\
\Gamma^n &:=& \left(\theta_2 + \sqrt{n^{1+ 2u}}  \cdot A \left((v-u) \tau_n , - \frac{\theta_{2}}{\sqrt{n}n^{v}} \right) \right) \cdot \sqrt{n}.
\end{eqnarray*}
{Using a similar argument as in the proof of  \eqref{eq:A-asy}, we find for $0<u<v \le T$,
\[
\lim_{n \rightarrow \infty} \hat \theta_1^{(n)} = \theta_1 + \theta_2, \quad \text{and} \quad
\lim_{n \rightarrow \infty}\Gamma^n = 0. \] }
It immediately follows that for $0<u<v \le T$,
\begin{equation*}
\lim_{n \rightarrow \infty} \mathbb{E}\left[e^{-\frac{\theta_{1}}{\sqrt{n}n^{u}}Z_{u\tau_n}
-\frac{\theta_{2}}{\sqrt{n}n^{v}}Z_{v\tau_n}}\bigg|Z_{0}=n\right] \cdot e^{\theta_{1}\sqrt{n}
+\theta_{2}\sqrt{n}} = \exp\left( \frac{1}{2} \frac{\alpha^{2}}{\alpha-\beta} {(\theta_1 + \theta_2)^2} \right),
\end{equation*}
which implies that
\begin{equation*}
\left(\frac{Z_{u\tau_n}-n^{1+u}}{\sqrt{n^{1+2u}}}, \frac{Z_{v\tau_n}-n^{1+v}}{\sqrt{n^{1+2v}}} \right)\rightarrow (Y(u), Y(v)),
\quad \text{weakly as $n \rightarrow \infty$},
\end{equation*}
where $Y(u) = Y(v) := \xi$ for $v>u>0$, and $\xi$ is a normal random variable with mean zero and variance $\frac{\alpha^{2}}{\alpha-\beta}$. It is clear that $Y$ is a continuous process
on $D[t,T]$ for any $0<t < T$. Thus by Proposition 1.2 in Aldous \cite{Aldous1989}, we deduce that
\begin{equation*}
\left\{\frac{Z_{s \tau_{n}}-n^{1+s}}{\sqrt{n^{1+2s}}}: s \in [t, T] \right\} \rightarrow Y,
\end{equation*}
weakly on $D[t,T]$.

We next prove \eqref{eq:Nlargetime-super-clt}, i.e., the functional central limit theorem for rescaled jump process $N$
in the super--critical case.
Recall from \eqref{eq:N-Z} that
\begin{equation*}
\frac{N_{s\tau_n}-\frac{n^{1+s}-n}{\alpha-\beta}}{n^{\frac{1}{2}+s}}
=\frac{Z_{s\tau_n}-n^{1+s}}{\alpha n^{\frac{1}{2}+s}}
+\frac{\beta}{\alpha}\tau_n\int_{0}^{s}n^{u-s}\frac{Z_{u\tau_n}-n^{1+u}}{n^{\frac{1}{2}+u}}du,
\end{equation*}
and that $Y_s \equiv \xi$ for $s>0$, where $\xi$ is a normal random variable with mean $0$ and variance $\frac{\alpha^{2}}{\alpha-\beta}$. In addition, note that
\begin{eqnarray*}
\tau_n\int_{0}^{s}n^{u-s} du = \frac{1}{\alpha - \beta} (1- n^{-s}).
\end{eqnarray*}
Therefore we deduce that
\begin{equation*}
\left\{\frac{N_{s\tau_n}-\frac{n^{1+s}-n}{\alpha-\beta}}{n^{\frac{1}{2}+s}}: s \in [t,T] \right\}
\rightarrow\frac{\xi}{\alpha}+\frac{\beta}{\alpha} \cdot \frac{1}{\alpha - \beta} \xi ,
\end{equation*}
weakly on $D[t, T]$. The weak limit is simply $\frac{1}{\alpha - \beta} \xi$.
\end{proof}

\subsubsection{Proof of super--critical case when $\mu>0$}
\begin{proof}
Zhu \cite{ZhuThesis} showed that when $\alpha>\beta$,
$\lim_{t\rightarrow\infty}\frac{1}{t}\log N_{t}^{(1)}=\alpha-\beta$.
Note that $Z_{t}^{(1)}\leq\alpha N_{t}^{(1)}$ and for any $0<T<\frac{1}{2}$,
\begin{equation*}
\sup_{0\leq s\leq T}\frac{Z^{(1)}_{s\tau_n}}{\sqrt{n^{1+2s}}}
\leq\frac{\alpha N^{(1)}_{T\tau_n}}{\sqrt{n}}\rightarrow 0,
\end{equation*}
almost surely as $n\rightarrow\infty$. Similarly, we have
\begin{equation*}
\sup_{0\leq s\leq T}\frac{N^{(1)}_{s\tau_n}}{\sqrt{n^{1+2s}}}
\leq\frac{ N^{(1)}_{T\tau_n}}{\sqrt{n}}\rightarrow 0,
\end{equation*}
almost surely as $n\rightarrow\infty$. Hence the results follow from the decompositions $Z= Z^{(0)} + Z^{(1)}$ and $N= N^{(0)} + N^{(1)}$ as given in Section~\ref{sec:mu-relation}.
\end{proof}

\subsubsection{Proof of sub--critical case when $\mu =0$}
\begin{proof}
We prove the convergence of finite dimensional distributions of the rescaled $Z$ processes. To this end, we first show for fixed $s<1$, the sequence of moment generating functions
\begin{equation*}
\mathbb{E}\left[e^{-\theta\frac{Z_{st_n}-n^{1-s} }{\sqrt{n^{1-s}}}}\bigg|Z_{0}=n\right]
\end{equation*}
converges when $n\rightarrow \infty$. 
We will show that for any $\theta\in\mathbb{R}$ and $s<1$, $\frac{|\theta|}{\sqrt{n^{1-s}}}<\theta_{c}(st_{n})$
for sufficiently large $n$,
where $\theta_{c}(\cdot)$ is defined in \eqref{thetacEqn} so that
for any fixed $\theta\in\mathbb{R}$, and for any sufficiently large $n$,
\begin{equation*}
\mathbb{E}\left[e^{-\theta\frac{Z_{st_n}-n^{1-s} }{\sqrt{n^{1-s}}}}\bigg|Z_{0}=n\right]
=\exp \left({ A\left(st_n, -\frac{\theta}{\sqrt{n^{1-s}}} \right) \cdot n+ \theta \sqrt{n^{1-s}} } \right).
\end{equation*}
To see this, we recall the definition of $\theta_{c}(t)$ in \eqref{thetacEqn}.
We note that in the sub--critical case $\alpha<\beta$, 
the function $-\beta A+e^{\alpha A}-1$ is $0$ at $A=0$, and it is convex in $A\geq 0$
and its derivative is negative at $A=0$ since $\alpha<\beta$.
Thus, there exists a unique positive value $A_{c}$ so that $-\beta A_{c}+e^{\alpha A_{c}}-1=0$
and $\theta_{c}(t)\rightarrow A_{c}$ as $t\rightarrow\infty$.
Thus, for any $\theta\in\mathbb{R}$ and $s<1$, we have $\frac{|\theta|}{\sqrt{n^{1-s}}}<\theta_{c}(st_{n})$
for sufficiently large $n$.
Next, we note that $e^{(\alpha-\beta) s t_n} = n^{-s}$ and the quantity $-\frac{\theta}{\sqrt{n^{1-s}}}$
goes to zero as $n \rightarrow \infty$. Using a similar argument as in the proof of  \eqref{eq:A-asy}, we find
\begin{equation*}
\lim_{n\rightarrow \infty} \exp \left({ A\left(st_n, -\frac{\theta}{\sqrt{n^{1-s}}} \right) \cdot n+ \theta \sqrt{n^{1-s}} } \right)
 = \exp\left( \frac{1}{2} \frac{\alpha^{2}}{\beta - \alpha} {\theta^2} \right).
\end{equation*}
Hence, we obtain
\begin{equation*}
{\lim_{n\rightarrow \infty}} \mathbb{E}\left[e^{-\theta\frac{Z_{st_n}-n^{1-s} }{\sqrt{n^{1-s}}}}\bigg|Z_{0}=n\right] = \exp\left( \frac{1}{2} \frac{\alpha^{2}}{\beta - \alpha} {\theta^2} \right),
\end{equation*}
which implies that for fixed $s \in (0,1)$,
\begin{equation*}
\frac{Z_{st_{n}}-n^{1-s}}{\sqrt{n^{1-s}}} \rightarrow R_s, \quad \text{weakly as $n \rightarrow \infty$}
\end{equation*}
where $R_s$ is a normal random variable with mean zero and variance $\frac{\alpha^{2}}{\beta - \alpha}$.

We proceed to consider the dimension of two.
The proof for the general case of finite dimensions follows from a similar argument and hence is omitted. Fix $\theta_1, \theta_2 \in\mathbb{R}$. 
For any $0<u<v<1$, and any sufficiently large $n$, the moment generating function of the random vector $\left(\frac{Z_{ut_{n}}-n^{1-u}}{\sqrt{n^{1 - u}}}, \frac{Z_{vt_{n}}-n^{1 - v}}{\sqrt{n^{1 - v}}} \right)$ is given by
\begin{equation} \label{eq:mgf-2d-rescaled1}
\mathbb{E}\left[e^{-\frac{\theta_{1}}{\sqrt{n^{1-u}}}Z_{ut_{n}}
-\frac{\theta_{2}}{\sqrt{{n^{1-v}}}}Z_{vt_{n}}}\bigg|Z_{0}=n\right] \cdot e^{\theta_{1}\sqrt{n^{1-u}}
+\theta_{2}\sqrt{n^{1-v}}}.
\end{equation}
By Cauchy-Schwarz inequality,
\begin{align*}
&\mathbb{E}\left[e^{-\frac{\theta_{1}}{\sqrt{n^{1-u}}}Z_{ut_{n}}
-\frac{\theta_{2}}{\sqrt{{n^{1-v}}}}Z_{vt_{n}}}\bigg|Z_{0}=n\right]
\\
&\leq
\left(\mathbb{E}\left[e^{-\frac{2\theta_{1}}{\sqrt{n^{1-u}}}Z_{ut_n}}\bigg|Z_{0}=n\right]\right)^{\frac{1}{2}}
\left(\mathbb{E}\left[e^{-\frac{2\theta_{2}}{\sqrt{n^{1-v}}}Z_{vt_n}}\bigg|Z_{0}=n\right]\right)^{\frac{1}{2}}
<\infty,
\end{align*}
for any sufficiently large $n$. 
It can be directly computed that
\begin{eqnarray*}
\lefteqn{\mathbb{E}\left[e^{-\frac{\theta_{1}}{\sqrt{n^{1-u}}}Z_{ut_{n}}
-\frac{\theta_{2}}{\sqrt{n^{1-v}}}Z_{vt_{n}}}\bigg|Z_{0}=n\right]} \\
& =& \mathbb{E}\left[e^{-\frac{\theta_{1}}{\sqrt{n^{1-u}}}Z_{ut_{n}}}
\cdot \left[ e^{ A((v-u) t_n , - \frac{\theta_{2}}{\sqrt{n^{1-v}}} ) \cdot Z_{u t_n} } \right]\Big|Z_{0}=n \right].
\end{eqnarray*}
After rewriting, the moment generating function in \eqref{eq:mgf-2d-rescaled1} becomes,
\begin{equation*}
\mathbb{E}\left[e^{- \hat \theta_1^{(n)} \frac{Z_{ut_n}-n^{1-u} }{\sqrt{n^{1-u}}} }\bigg|Z_{0}=n\right] \cdot \exp(\Gamma^n)
\end{equation*}
where, with slight abuse of notations, 
\begin{eqnarray*}
\hat \theta_1^{(n)} &:=&\theta_1 - \sqrt{n^{1- u}}  \cdot A \left((v-u) t_n , - \frac{\theta_{2}}{\sqrt{n^{1-v}}} \right),\\
\Gamma^n &:=& {n^{1- u}}  \cdot A \left((v-u) t_n , - \frac{\theta_{2}}{\sqrt{n^{1-v}}} \right) + \theta_2 \cdot\sqrt{n^{1-v}} .
\end{eqnarray*}
{Using a similar argument as in the proof of  \eqref{eq:A-asy}, we find}
\[
\lim_{n \rightarrow \infty} \hat \theta_1^{(n)} = \theta_1 \quad \text{and} \quad
\lim_{n \rightarrow \infty}\Gamma^n = \frac{1}{2} \frac{\alpha^{2}}{\beta - \alpha} {\theta_2^2} .
\]
It immediately follows that
\begin{equation*}
\lim_{n \rightarrow \infty}\mathbb{E}\left[e^{-\frac{\theta_{1}}{\sqrt{n^{1-u}}}Z_{ut_{n}}
-\frac{\theta_{2}}{\sqrt{{n^{1-v}}}}Z_{vt_{n}}}\bigg|Z_{0}=n\right] \cdot e^{\theta_{1}\sqrt{n^{1-u}}
+\theta_{2}\sqrt{n^{1-v}}} = \exp\left( \frac{1}{2} \frac{\alpha^{2}}{\beta -\alpha} {(\theta_1^2 + \theta_2^2)} \right),
\end{equation*}
which implies that
\begin{equation*}
\left(\frac{Z_{ut_{n}}-n^{1-u}}{\sqrt{n^{1 - u}}}, \frac{Z_{vt_{n}}-n^{1 - v}}{\sqrt{n^{1 - v}}} \right)\rightarrow (R_u, R_v),
\quad \text{weakly as $n \rightarrow \infty$},
\end{equation*}
where $R_u$ and $R_v$ are independent normal random variables, both with mean zero and variance $\frac{\alpha^{2}}{\beta -\alpha}$. The proof is complete.
\end{proof}

\subsubsection{Proof of sub--critical case $\mu>0$}
\begin{proof}
Note that for $\beta>\alpha>0$, $\frac{N_{t}^{(1)}}{t}\rightarrow\frac{\mu}{1-\frac{\alpha}{\beta}}$ almost surely
as $t\rightarrow\infty$. Therefore, for any $0<T<1$,
\begin{equation*}
\sup_{0\leq s\leq T}\frac{Z_{st_{n}}^{(1)}}{\sqrt{n^{1-s}}}
\leq\frac{\alpha N_{Tt_{n}}^{(1)}}{\sqrt{n^{1-T}}}\rightarrow 0,
\end{equation*}
almost surely as $n\rightarrow\infty$. Hence the result follows from the decomposition $Z=Z^{(0)}+Z^{(1)}$ as described in Section~\ref{sec:mu-relation}.
\end{proof}


\section{Appendix}
\subsection{Proof of Proposition~\ref{prop:moments-noncritial}}

\begin{proof}
The process $Z$ is a piecewise deterministic
Markov process as defined in Davis \cite{Davis1984}; see also Chapter 11 in Rolski et al. \cite{Rolski2009}.
Suppose $f$ is any continuous function with at most polynomial growth on $(0, \infty)$, i.e., there exists some positive integer $k$ such that $|f(x)| \le x^k$ for all $x \in (0, \infty)$. Recall that {$dZ_{t}=-\beta Z_{t}dt+\alpha dN_{t}$}. Define
\begin{equation} \label{eq:generator}
\mathcal{A}f(z)=-\beta z\frac{\partial f}{\partial z}
+z[f(z+\alpha)-f(z)].
\end{equation}
We first verify that $f$ is in the domain of
the full generator of the Markov process $\{Z_t : t \ge 0\} $\footnote{The domain of the infinitesimal generator of a Markov process is always contained in the domain of its full generator.}. That is,
\begin{equation*}
\left\{f(Z_{t}) - f(Z_{0}) - \int_{0}^{t}\mathcal{A}f(Z_{s})ds: t \ge 0 \right\}
\end{equation*}
is a martingale with respect to the natural filtration generated by $Z$. (See Section 11.1.4 in Rolski et al. (2009) for the
definition of full generator and further details). We use Theorem 11.2.2 in Rolski et al. \cite{Rolski2009} and check the three conditions there. Since the boundary
set of the piecewise deterministic Markov process $Z$ is empty ($Z$ never hits zero), and the sample path of $Z$ is absolutely continuous between jumps, it suffices to check that for each $t \ge 0$,
\begin{eqnarray*}
\mathbb{E}_{Z_0} \left(\sum_{\tau_i \le t } \left| f(Z_{\tau_i}) - f(Z_{\tau_i -}) \right| \right) < \infty,
\end{eqnarray*}
where $\tau_{i}$'s are the jump epochs of the process $Z$ andand $ \mathbb{E}_{Z_0} [\cdot] := \mathbb{E} [\cdot |Z_0]$ for given $Z_0$. Direct computation yields
\begin{eqnarray*}
\mathbb{E}_{Z_0} \left(\sum_{\tau_i \le t } \left| f(Z_{\tau_i}) - f(Z_{\tau_i -}) \right| \right) &=& \mathbb{E}_{Z_0}
\left[\int_{0}^t |f(Z_{s-}+\alpha)-f(Z_{s-})| dN_s\right]\\
&\le& 2\mathbb{E}_{Z_0} \left[\int_{0}^t (Z_{s-}+\alpha)^k dN_s\right]\\
&\le& 2\mathbb{E}_{Z_0} \left[\int_{0}^t (Z_{0}+\alpha N_{s-}+\alpha)^k dN_s\right]\\
&\le& 2^{k} \mathbb{E}_{Z_0} \left[\int_{0}^t ((Z_0+\alpha)^{k} + \alpha^{k} N_{s-}^{k}) dN_{s}\right] \\
&\le& 2^{k} \left((Z_{0}+\alpha)^{k} \mathbb{E}_{Z_0} [N_{t}] + \alpha^{k}\mathbb{E}_{Z_0} [N_{t}^{k+1}] \right)< \infty,
\end{eqnarray*}
where we have used the facts that $f$ is of at most polynomial growth, $Z_t \le Z_0 + \alpha N_t$ for all $t$, and the moments of $N_t$ are finite (see, e.g., \cite{ZhuCLT})
\footnote{Note that in \cite{ZhuCLT}, it is proved
that $\mathbb{E}[e^{\theta N_{t}}]<\infty$ for any sufficiently small $\theta>0$, 
for $\mathbb{E}$ being the expectation for the stationary Hawkes process.
In our setting, by tower property, $\mathbb{E}[e^{\theta N_{t}}]=\mathbb{E}[\mathbb{E}_{Z_{0}}[e^{\theta N_{t}}]]<\infty$,
which implies that $\mathbb{E}_{Z_{0}}[e^{\theta N_{t}}]<\infty$ for a.e. $Z_{0}$. 
Since $\mathbb{E}_{Z_{0}}[e^{\theta N_{t}}]$ is monotonic in $Z_{0}$, 
we conclude that $\mathbb{E}_{Z_{0}}[e^{\theta N_{t}}]<\infty$ for every $Z_0$.
Hence, it follows that $\mathbb{E}_{Z_{0}}[N_{t}^{k}]<\infty$ for every $k\in\mathbb{N}$.}. The martingale property then implies
\begin{equation}\label{eq:dynkin}
\mathbb{E}_{Z_0}[f(Z_{t})]=f(Z_{0})+\mathbb{E}_{Z_0}\left[\int_{0}^{t}\mathcal{A}f(Z_{s})ds\right],
\end{equation}
where $\mathcal{A}f$ is given in \eqref{eq:generator}.

To compute the first moment, we apply $f(z)=z$ in \eqref{eq:dynkin} and obtain
\begin{equation*}
\mathbb{E}_{Z_0}[Z_{t}]=Z_{0} + \mathbb{E}_{Z_0} \left[ \int_{0}^{t} (\alpha - \beta) Z_s ds \right] = Z_{0} + \int_{0}^{t} (\alpha - \beta) \cdot \mathbb{E}_{Z_0}[Z_{s}] ds,
\end{equation*}
which implies
\begin{equation*}
\frac{d}{dt} \mathbb{E}_{Z_0}[Z_{t}] = (\alpha - \beta)\cdot\mathbb{E}_{Z_0}[Z_{t}].
\end{equation*}
Thus we have
\begin{equation*}
\mathbb{E}_{Z_0}[Z_{t}]=Z_{0}e^{(\alpha-\beta)t}.
\end{equation*}

Now we compute second moments. Applying $f(z)=z^2$ in \eqref{eq:dynkin}, we find
\begin{align*}
\mathbb{E}_{Z_0}[Z_{t}^{2}]&=Z_{0}^{2}+2(\alpha-\beta)\int_{0}^{t}\mathbb{E}_{Z_0} [Z_{s}^{2}]ds
+\alpha^{2}\int_{0}^{t}\mathbb{E}_{Z_0} [Z_{s}]ds
\\
&=Z_{0}^{2}+2(\alpha-\beta)\int_{0}^{t}\mathbb{E}_{Z_0} [Z_{s}^{2}]ds
+\alpha^{2}Z_{0}\int_{0}^{t}e^{(\alpha-\beta)s}ds.
\nonumber
\end{align*}
Solving this equation yields
\begin{align*}
\mathbb{E}_{Z_0} [Z_{t}^{2}] =
\begin{cases}
Z_{0}^{2}+\alpha^{2}Z_{0}t , & \qquad \text{if} \quad \alpha = \beta,\\
Z_{0}^{2}e^{2(\alpha-\beta)t}
+\frac{\alpha^{2}Z_{0}}{\alpha-\beta}(e^{2(\alpha-\beta)t}-e^{(\alpha-\beta)t}),      & \qquad \text{if} \quad \alpha \ne \beta.
\end{cases}
\end{align*}

Next we compute the third moments. We apply $f(z)=z^3$ in \eqref{eq:dynkin}. If $\alpha = \beta$, we immediately get
\begin{align*}
\mathbb{E}_{Z_0} [Z_{t}^{3}]
&=Z_{0}^{3}+3\alpha^{2}\int_{0}^{t}\mathbb{E}_{Z_0} [Z_{s}^{2}]ds
+\alpha^{3}\int_{0}^{t}\mathbb{E}_{Z_0} [Z_{s}]ds
\\
&=Z_{0}^{3}+3\alpha^{2}Z_{0}^{2}t+\frac{3}{2}\alpha^{4}Z_{0}t^{2}+\alpha^{3}Z_{0}t.
\nonumber
\end{align*}
When $\alpha \ne \beta$, we have
\begin{align*}
\mathbb{E}_{Z_0}[Z_{t}^{3}]&=Z_{0}^{3}+3(\alpha-\beta)\int_{0}^{t}\mathbb{E}_{Z_0}[Z_{s}^{3}]ds
+3\alpha^{2}\int_{0}^{t}\mathbb{E}_{Z_0}[Z_{s}^{2}]ds
+\alpha^{3}\int_{0}^{t}\mathbb{E}_{Z_0}[Z_{s}]ds
\\
&=Z_{0}^{3}+3(\alpha-\beta)\int_{0}^{t}\mathbb{E}_{Z_0}[Z_{s}^{3}]ds
+3\alpha^{2}Z_{0}^{2}\int_{0}^{t}e^{2(\alpha-\beta)s}ds
\nonumber
\\
&\qquad
+3\alpha^{2}\frac{\alpha^{2}Z_{0}}{\alpha-\beta}\int_{0}^{t}(e^{2(\alpha-\beta)s}-e^{(\alpha-\beta)s})ds
+\alpha^{3}Z_{0}\int_{0}^{t}e^{(\alpha-\beta)s}ds
\nonumber
\end{align*}
Therefore,
\begin{equation*}
\frac{d}{dt}\mathbb{E}_{Z_0}[Z_{t}^{3}]
=3(\alpha-\beta)\mathbb{E}_{Z_0}[Z_{t}^{3}]
+\left(3\alpha^{2}Z_{0}^{2}+\frac{3\alpha^{4}Z_{0}}{\alpha-\beta}\right)e^{2(\alpha-\beta)t}
+\left(\alpha^{3}Z_{0}-\frac{3\alpha^{4}Z_{0}}{\alpha-\beta}\right)e^{(\alpha-\beta)t},
\end{equation*}
which yields that
\begin{align*}
\mathbb{E}_{Z_0}[Z_{t}^{3}]e^{-3(\alpha-\beta)t}-Z_{0}^{3}
&=\left(\frac{3\alpha^{2}Z_{0}^{2}}{\alpha-\beta}+\frac{3\alpha^{4}Z_{0}}{(\alpha-\beta)^{2}}\right)
(1-e^{-(\alpha-\beta)t})
\\
&\qquad\qquad
+\left(\frac{\alpha^{3}Z_{0}}{2(\alpha-\beta)}-\frac{3\alpha^{4}Z_{0}}{2(\alpha-\beta)^{2}}\right)(1-e^{-2(\alpha-\beta)t}),
\nonumber
\end{align*}
which implies that
\begin{align*}
\mathbb{E}_{Z_0}[Z_{t}^{3}]
&=\left(Z_{0}^{3}+\frac{3\alpha^{2}Z_{0}^{2}}{\alpha-\beta}
+\frac{\alpha^{3}Z_{0}}{2(\alpha-\beta)}+\frac{3\alpha^{4}Z_{0}}{2(\alpha-\beta)^{2}} \right)
e^{3(\alpha-\beta)t}
\\
&\qquad
-\left(\frac{3\alpha^{2}Z_{0}^{2}}{\alpha-\beta}+\frac{3\alpha^{4}Z_{0}}{(\alpha-\beta)^{2}}\right)
e^{2(\alpha-\beta)t}
\nonumber
\\
&\qquad\qquad
-\left(\frac{\alpha^{3}Z_{0}}{2(\alpha-\beta)}-\frac{3\alpha^{4}Z_{0}}{2(\alpha-\beta)^{2}}\right)
e^{(\alpha-\beta)t}.
\nonumber
\end{align*}
The proof is therefore complete.
\end{proof}

\subsection{Proof of Proposition~\ref{prop:ineq}}
\begin{proof}
First, recall from \eqref{MartingaleRep} that
\begin{equation}\label{eq:tmp}
Z_{t}-ne^{(\alpha-\beta)t}=e^{(\alpha-\beta)t}\alpha\int_{0}^{t}e^{-(\alpha-\beta)s}dM_{s}.
\end{equation}
We note that $\{\int_{0}^{t}e^{-(\alpha-\beta)s}dM_{s}: t \in [0, T]\}$ is a martingale
with predictable quadratic variation $\int_{0}^{t}e^{-2(\alpha-\beta)s}d\langle M\rangle_{s}$,
where $\langle M\rangle_{t}=\int_{0}^{t}Z_{s}ds$ is the predictable quadratic variation
of the martingale $M$.
By Burkholder-Davis-Gundy inequality, we have
\begin{eqnarray*}
\mathbb{E}\left[\left(\int_{0}^{t}e^{-(\alpha-\beta)s}dM_{s}\right)^{4}\right]
&\leq&\mathbb{E}\left[\left(\sup_{0\leq s\leq t}\left|\int_{0}^{s}e^{-(\alpha-\beta)u}dM_{u}\right|\right)^{4}\right]\\
&\leq& C\cdot\mathbb{E}\left[\left(\int_{0}^{t}e^{-2(\alpha-\beta)s}d\langle M\rangle_{s}\right)^{2}\right],
\end{eqnarray*}
where $C$ is a positive constant.
Hence, when $\alpha\neq\beta$, we deduce from Proposition~\ref{prop:moments-noncritial} that
\begin{align}
&\mathbb{E}\left[\left(\int_{0}^{t}e^{-(\alpha-\beta)s}dM_{s}\right)^{4}\right]
\label{eq:mart1} \\
&\leq C\cdot\mathbb{E}\left[\left(\int_{0}^{t}e^{-2(\alpha-\beta)s}Z_{s}ds\right)^{2}\right]
\nonumber(\text{Burkholder-Davis-Gundy inequality})
\\
&\leq Ct\mathbb{E}\left[\int_{0}^{t}e^{-4(\alpha-\beta)s}Z_{s}^{2}ds\right]
\nonumber(\text{Cauchy-Schwarz inequality})
\\
&=Ct\int_{0}^{t}e^{-4(\alpha-\beta)s}\mathbb{E}[Z_{s}^{2}]ds
\nonumber
\\
&=Ct\int_{0}^{t}e^{-4(\alpha-\beta)s}\left[Z_{0}^{2}e^{2(\alpha-\beta)s}
+\frac{\alpha^{2}Z_{0}}{\alpha-\beta}(e^{2(\alpha-\beta)s}-e^{(\alpha-\beta)s})\right]ds
\nonumber
\\
&=Ct\left[\frac{n^{2}}{2(\alpha-\beta)}(1-e^{-2(\alpha-\beta)t})
+\frac{\alpha^{2}n}{2(\alpha-\beta)^{2}}(1-e^{-2(\alpha-\beta)t})
-\frac{\alpha^{2}n}{3(\alpha-\beta)^{2}}(1-e^{-3(\alpha-\beta)t})\right].
\nonumber
\end{align}
Hence we have established \eqref{ineq:mart}.
Similarly, when $\alpha = \beta$, we can obtain that 
\begin{align}
\mathbb{E}\left[\left(\int_{0}^{t}e^{-(\alpha-\beta)s}dM_{s}\right)^{4}\right]
& \le Ct\int_{0}^{t}e^{-4(\alpha-\beta)s}\mathbb{E}[Z_{s}^{2}]ds
\nonumber
\\
&=Ct \left(Z_0^2 t + \frac{1}{2} \alpha^2 t^2 Z_0\right)=Ct \left(n^2 t + \frac{1}{2} \alpha^2 t^2 n\right).
\label{eq:mart2}
\end{align}
Using \eqref{eq:tmp}, \eqref{eq:mart1} and \eqref{eq:mart2}, we infer that the inequality \eqref{UpperBoundI} also holds. 

To prove \eqref{eq:mart-upperbound}, we can use a similar argument as above. When $\alpha \ne \beta,$ one can readily verify that
\begin{align*}
&\sup_{\delta \le t \le T}\mathbb{E}\left[\left(\int_{t-\delta}^{t}e^{-(\alpha-\beta)s}dM_{s}\right)^{4}\right]\\
&\le C \delta \sup_{\delta \le t \le T} \int_{t -\delta}^{t}e^{-4(\alpha-\beta)s}\left[Z_{0}^{2}e^{2(\alpha-\beta)s}
+\frac{\alpha^{2}Z_{0}}{\alpha-\beta}\left(e^{2(\alpha-\beta)s}-e^{(\alpha-\beta)s}\right)\right]ds\\
&\leq C n^{2}\delta^2.
\end{align*}
The proof of the case $\alpha = \beta$ is similar and is thus omitted. The proof is therefore complete.
\end{proof}

\section*{Acknowledgements}
We are very grateful to the AE and the referee for their careful reading of the manuscript
and very helpful comments and suggestions.
Xuefeng Gao acknowledges support from Hong Kong RGC ECS grant 24207015 
and CUHK Direct Grants for Research with project codes 4055035 and 4055054.
Lingjiong Zhu is grateful to the support from NSF Grant DMS-1613164.
We would like to thank Peter Carr, Rama Cont, Jasmine Foo, Fuqing Gao, 
Michel Mandjes and Frederi G. Viens for helpful discussions.


\end{document}